\documentclass[10pt]{article}
\usepackage[T1]{fontenc}
\usepackage{euscript}
\usepackage{amsmath}
\usepackage{amssymb}
\usepackage{latexsym}
\usepackage{fancyhdr}

\usepackage{diagrams}
\include{amsfonts}

\newarrow{equal}=====
\newarrow{to}---->
\newarrow{onto}----{>>}
\newarrow{embed}>--->
\newarrow{TeXonto}{-}{-}{-}{-}{->>}
\newarrow{TeXto}----{->}
\newarrow{TeXembed}{>}---{->}
\newarrow{TeXembedd}{>}---{->>}

\newtheorem{df}{Definition}[section]      
\newtheorem{thm}[df]{Theorem}             
\newtheorem{prop}[df]{Proposition}
\newtheorem{cor}[df]{Corollary}

\newtheorem{lem}[df]{Lemma}

\newtheorem{rmk}[df]{Remark}

\newcommand{\pf}{\noindent{\sc Proof.}\ }

\newcommand{\boom}{\quad\lower3pt\hbox{\vrule height1.1ex width .9ex depth -.2ex}
                    \vskip9pt}

\renewcommand{\mathcal}[1]{\EuScript{#1}}

\let\phi=\varphi

\let\da=\tau

\newcommand{\ad}{\mathop{\rm ad}\nolimits}
\newcommand{\Ad}{\mathop{\rm Ad}\nolimits}

\def\Aut{\mathop{\rm Aut}}
\newcommand{\coker}{\mathop{\rm coker}}
\newcommand{\CDO}{\mathop{\rm CDO}}
\def\Der{\mathop{\rm Der}}

\def\im{\mathop{\rm im}}
\def\Int{\mathop{\rm Int}}
\def\Inn{\mathop{\rm Inn}}

\def\Opext{\mathop{{\cal O}{\rm pext}}}
\def\Obs{\mathop{\rm Obs}}

\def\OutDO{\mathop{\rm OutDO}}

\def\xm{\mathop{\rm xm}}

\def\R{\mathbb{R}}    
\def\chigh{{\raise1.5pt\hbox{$\chi$}}}

\def\actgpd{\mathbin{\hbox{$<\kern-.4em\mapstochar\kern.4em$}}}
\def\ractgpd{\mathbin{\hbox{$\mapstochar\kern-.3em>$}}}

\newcommand{\gog}{\mathfrak{g}}
\newcommand{\hoh}{\mathfrak{h}}

\newcommand{\co}{\colon\thinspace}                

\def\act{\mathbin{\hbox{$<\kern-.4em\mapstochar\kern.4em$}}}
\def\ract{\mathbin{\hbox{$\mapstochar\kern-.3em>$}}}

\let\Bar=\overline
\let\Hat=\widehat
\let\Tilde=\widetilde

\renewcommand{\to}{\longrightarrow}

\def\gpd{\,\lower1pt\hbox{$\longrightarrow$}\hskip-.24in\raise2pt
             \hbox{$\longrightarrow$}\,}

\newcommand{\surj}{-\!\!\!-\!\!\!-\!\!\!\gg}
\newcommand{\inj}{>\!\!\!-\!\!\!-\!\!\!-\!\!\!>}

\begin{document}

\title{{\bf CROSSED MODULES AND THE INTEGRABILITY OF LIE BRACKETS}
} 

\author{Iakovos Androulidakis\\
        Institut f\'{u}r Mathematik\\
        Universit\"{a}t Z\"{u}rich\\
        190 Winterthurerstrasse\\
        Zurich CH-8057\\
        Switzerland\\
        {\sf iakovos.androulidakis@math.unizh.ch}}

\date{{\sf \LaTeX version of \today}}

\maketitle

\begin{abstract}
We show that the integrability obstruction of a transitive Lie algebroid coincides with the lifting obstruction of a crossed module of groupoids associated naturally with the given algebroid. Then we extend this result to general extensions of integrable transitive Lie algebroids by Lie algebra bundles. Such a lifting obstruction is directly related with the classification of extensions of transitive Lie groupoids. We also give a classification of such extensions which differentiates to the classification of transitive Lie algebroids discussed in \cite{KCHM:new}.
\end{abstract}

\tableofcontents


\section*{Introduction}

A Lie algebroid may be thought of as a generalisation of the tangent bundle of a manifold. More formally, it is a vector bundle $A \to M$ whose module of sections is equipped with a Lie bracket preserved by a certain vector bundle morphism $q : A \to TM$ (anchor). Its global counterpart is a Lie groupoid, which, roughly speaking, is a small category such that every arrow is invertible, with a suitable smooth structure. Both notions play an important role in a number of different areas of geometry: In differential geometry the notions of connection and holonomy are naturally formulated in terms of (transitive) algebroids and groupoids; in Poisson geometry, there is a loose duality between Lie algebroids and Poisson manifolds; in noncommutative geometry, groupoids provide certain $C^*$-algebras which substitute the often pathological leaf space of a foliation.

Any Lie algebroid over a point is a Lie algebra and every Lie groupoid whose space of units is a point is a Lie group. Groupoids and algebroids therefore generalise the classical Lie groups-Lie algebras apparatus. We cannot, however,  generalise to the "-oids" context Lie's Third Theorem, namely not every Lie algebroid integrates to a Lie groupoid. Explicit obstructions were given by Crainic and Fernandes in \cite{Crainic-Fernandes}. The integrability problem comes up quite often, in problems arising in both mathematics and physics. For instance, if a Lie algebroid $A$ integrates to a Lie groupoid $G$ then the commutators of $C^*(G)$ approximate asymptotically the natural Poisson bracket on the functions on $A^*$. This is a deformation quantization of this bracket. In another direction the fiberwise additive groupoid structure of  $A$ may be compined with $G$ to the tangent groupoid $A \times \{ 0 \} \cup G \times \R^*$ over $M \times \R$ to give rise to the analytic index map along the leaves of the foliation defined by the image of the anchor map (see \cite{Connes}).

This paper is concerned with the integrability problem of the class of {\it transitive} Lie algebroids, namely the ones whose anchor map is surjective, therefore the foliation they induce on $M$ has only one leaf. Many important Lie algebroids are transitive, such as the one associated naturally with a symplectic manifold and the Atiyah sequence of a principal bundle. Also, the restriction of any algebroid $A \to M$ on a leaf of the associated foliation is transitive. In the recent preprint \cite{Li} L.-C. Li showed that Hamiltonian systems may be realised as certain coboundary dynamical Poisson groupoids which are transitive.

We would like to draw special attention to the integrability problem for the algebroid $A_{\omega}$ of a symplectic manifold $(M,\omega)$. This problem is equivalent to the integrality problem for $\omega$ (see \cite[II\S8.1]{KCHM:new}). Conceptually, $\omega$ is integral iff it can be realised as the curvature of a connection on a principal bundle $P(M,G)$. If such data exists then the transitive groupoid $\frac{P \times P}{G} \gpd M$ corresponding to $P(M,G)$ integrates the Lie algebroid $A_{\omega}$  From this point of view geometric prequantization becomes an instance in the integrability of transitive Lie algebroids. But geometric prequantization is a two-step process, namely one needs to
\begin{enumerate}
\item Examine whether the symplectic form $\omega$ is integral.
\item If so, classify all principal bundles $P(M,G)$ with connection whose curvature is $\omega$.
\end{enumerate}

Transitive Lie groupoids and Lie algebroids were studied thoroughly by Mackenzie in the 1980's, with his work culminating to \cite[I, II]{KCHM:new}. The material there is influenced by the following observation: In the transitive case, since the anchor map is surjective the bundle $L = \ker q$ is a Lie algebra bundle . Therefore any transitive Lie algebroid is in fact an extension
\begin{eqnarray}\label{extn:1}
L \inj A \surj TM
\end{eqnarray}
When $M$ is simply connected, the integrability obstruction is a certain element of $\check{H}^{2}(M,\Tilde{G})$, where $\Tilde{G}$ is the connected and simply connected Lie group integrating the fiber of $L$. On the global level, every transitive Lie groupoid $\Omega \gpd M$ can be written as an extension
\begin{eqnarray}\label{extn:2}
F \inj \Omega \stackrel{(t,s)}{\surj} M \times M
\end{eqnarray}
of the pair groupoid $M \times M \gpd M$ by the Lie group bundle $F \to M$ of isotropies. In terms of prequantization, the integrability obstruction amounts to $(i)$ above. The analogous classification of $(ii)$ for transitive Lie groupoids was given in \cite{Mackenzie:1989} using the notion of a {\em crossed module}.

Crossed modules of Lie groups were introduced by Whitehead \cite{Whitehead:1949} in the context of homotopy theory. They were used later in the classification of extensions of Lie groups to overcome the problem that the group of inner automorphisms of a Lie group may not be closed in the full automorphism group. Namely, consider an extension of groups 
\[
N \stackrel{\iota}{\inj} H \stackrel{\pi}{\surj} G.
\]
If $N$ is abelian, the map $\rho : G \rightarrow Aut(N)$ defined by $\rho(g) = I_{h}\mid_{N}$, where $h$ is any element of $H$ such that $\pi(h) = g$, is a well defined representation. Now if $N$ is non-abelian, the automorphism $\rho$ is no longer well defined. The usual way around this problem is to consider the map $\rho : G \rightarrow Out(N) = \frac{Aut(N)}{Inn(N)}$, given by $g \mapsto \langle I_{h}\mid_{N} \rangle$, where $\pi(h) = g$. Here $Inn(N)$ is the group of inner automorphisms of $N$. This is a well defined morphism, called the {\em abstract kernel} of the original extension, and there is a standard classification of such extensions with a prescribed abstract kernel. When one is dealing with Lie groups though, the previous approach is problematic, because $Inn(N)$ need not be closed in $Aut(N)$, and the smoothness of the representation $\rho : G \rightarrow Out(N)$ has no longer a meaning. An alternative approach, circumventing this problem, was given by Mackenzie in \cite{Mackenzie:1989}, using crossed modules of Lie groups.
The link of crossed modules with cohomology can be traced back to \cite{Huebschmann}. Crossed modules of Lie groupoids were considered by Brown and Spencer \cite{Brown-Spencer:1976}, Brown and Higgins \cite{Brown-Higgins:1981}, and by Mackenzie in \cite{Mackenzie:1989} to classify principal bundles with prescribed gauge group bundle. In another direction, crossed modules are arise in string theory.

Here we introduce the notion of a crossed module of Lie algebroids and discuss its relevance with the integrability of transitive Lie algebroids. The main results of this paper are:
\begin{itemize}
\item For any transitive Lie algebroid the integrability obstruction coincides with the lifting obstruction of a certain crossed module of groupoids naturally associated with the given Lie algebroid.
\item For extensions of integrable transitive Lie algebroids by Lie algebra bundles, the integrability obstruction coincides with the lifting obstruction of a natural crossed module of groupoids naturally associated with the given extension. 
\item In both cases we classify cohomologically the possible lifts in case the obstruction vanishes.
\end{itemize}
This way we get a unified approach to the integrability problem. Namely by reformulating it to a lifting problem in the language of crossed modules we deal with both the integrability obstruction and the classification of the integrating groupoids. This is made clearer when one considers the second result above 
which concerns extensions of the form
\begin{eqnarray}\label{extn:3}
L \to A \to A\Omega
\end{eqnarray}
The integrability obstruction of extensions (\ref{extn:3}) was studied by Mackenzie in \cite{Mackenzie:1989IOELA}. The problem was reduced to the case of extensions of the form (\ref{extn:1}) plus the action of a Lie group. It was shown that such extensions over $M$ are equivalent to Lie algebroids $A' \to P$, where $P$ is the total space of a principal bundle $P(M,G)$ (the one corresponding to the groupoid $\Omega$), together with an action of the group $G$ by automorphisms. Such structures were called PBG-algebroids. Here we adopt this equivalence. The integrability obstruction of extensions was given in \cite{Mackenzie:1989IOELA}, but not the classification of the possible integrating groupoids. The problem with generalising the approach of \cite{Mackenzie:1989} is that PBG-groupoids involve a group action, not however recorded by the usual \v{C}ech cohomology used in \cite{Mackenzie:1989}. This difficulty was overcome by the introduction of {\em isometablic} cohomology in \cite{Androulidakis:solo1} and \cite{Androulidakis:2004CELG}. In this paper we reformulate isometablic cohomology to a form that classifies extensions (rather than principal bundles with an extra group action) in order to classify appropriately the lifts of crossed modules of PBG-groupoids. The use of crossed modules forces one to think of such lifts as extensions. Namely we prove that:
\begin{itemize}
\item Extensions of transitive Lie groupoids $F \inj \Omega \surj \Omega'$ are classified by suitably equivariant pairs $(\Tilde{\chi},\Tilde{\alpha})$ where $\Tilde{\alpha}_{ij} : P_{ij} \to \Aut{H}$ is the cocycle classifying to $F$ and $\chi_{ij} : P_{ij} \times P_{ij} \to H$ is a morphism of Lie groupoids, such that
\begin{enumerate}
\item $\Tilde{\chi}_{ik}(u,v) = \Tilde{\chi}_{ij}(u,v) \cdot \Tilde{\alpha}_{ij}(u)(\tilde{\chi}_{jk}(u,v))$
\item $\Tilde{\alpha}_{ij}(u) = I_{\Tilde{\chi}_{ij}(u,\cdot)}$
\end{enumerate}
\end{itemize}

Such a classification settles yet another matter. In \cite{KCHM:new} transitive Lie algebroids were viewed as extensions (\ref{extn:1}) and classified by pairs $(\chi,\alpha)$, where $\alpha_{ij} : U_{ij} \to \Aut(\gog)$ is the cocycle corresponding to the Lie algebra bundle $L$ and $\chi_{ij} : TU_{ij} \times TU_{ij} \to U_{ij} \times \gog$ are Maurer-Cartan 2-forms satisfying
\begin{enumerate}
\item $\chi_{ik} = \chi_{ij} + \alpha_{ij}(\chi_{jk})$
\item $\Delta(\alpha_{ij}) = \ad \circ \chi_{ij}$
\end{enumerate}
(here $\Delta$ stands for the Darboux derivative). On the other hand, a transitive Lie groupoid $\Omega \gpd M$ was classified by the \v{C}ech cocycle $\{ s_{ij} \}$ which classifies its corresponding principal bundle. It is not clear how the cocycle  $\{ s_{ij} \}$ differentiates to the data $(\chi,\alpha)$.





The paper is structured as follows: Section 1 is a short overview of crossed modules of Lie groupoids. In section 2 we introduce crossed modules of Lie algebroids and show that they generalise the notion of coupling discussed in \cite{KCHM:new}. We also discuss the lifting problem. Section 3 gives our method to recover the integrability obstruction using crossed modules. In section 4 we give the equivariant version of crossed modules, on the algebroid and the groupoid level, and discuss the integration/differentiation process. We also show that it suffices to consider only the simplest case of such crossed modules. Section 5 gives the classification of PBG-groupoids we discussed above. Sections 6 and 7 give the lifting obstruction and classify the lifts in case the obstruction vanishes.

As far as the non-transitive case is concerned, a classification of extensions of Lie groupoids which induce a regular foliation on the base manifold was given by Moerdijk \cite{Moerdijk:classification}. Our method to produce a crossed module of Lie groupoids from a given transitive Lie algebroid may easily be generalised in the non-transitive but regular case. In section 3 we discuss this and some implications in quantization. It would be very interesting to investigate whether a crossed modules approach can be applied to more general regular extensions.

\section*{Acknowledgments}

The idea of linking the integrability obstruction with crossed modules was developed during the author's PhD studies, in private communication with K. Mackenzie. We would like to thank the referee for various comments concerning the presentation of this material.

\section{Crossed modules of Lie groupoids}

This section and the next intend to set up the notion of crossed module for Lie groupoids and for Lie algebroids. We start here at the groupoid level, briefly recalling and clarifying material from \cite{Mackenzie:1989}. 

\begin{df} \label{def:gpdrep}
Let $\Omega$ be a Lie groupoid on a base manifold $M$ and let $(F,\pi,M)$ be a Lie group bundle
on $M$. A {\em representation} of $\Omega$ on $F$ is a smooth map $\rho : \Omega * F
\rightarrow F$, where $\Omega * F$ is the pullback manifold $\{ (\xi,f) \in \Omega \times F
: \alpha(\xi) = \pi(f) \}$, such that
\begin{enumerate}
\item $\pi(\rho(\xi,f)) = \beta(\xi)$ for $(\xi,f) \in \Omega * F$;
\item $\rho(\eta,\rho(\xi,f)) = \rho(\eta \xi,f)$ for all $f,\eta,\xi$ such that $(\xi,f) \in
\Omega * F$ and $(\eta,\xi) \in \Omega * \Omega$;
\item $\rho(1_{\pi(f)},f) = f$ for all $f \in F$;
\item $\rho(\xi) : F_{\alpha(\xi)} \rightarrow F_{\beta(\xi)},~f \mapsto \rho(\xi,f)$ is a
Lie group isomorphism for all $\xi \in \Omega$.
\end{enumerate}
\end{df}

Representations of groupoids on fibered manifolds were introduced by Ehresmann. One can also 
think of a groupoid representation as a Lie groupoid morphism $\Omega \rightarrow \Phi(F)$, where 
$\Phi(F)$ is Lie groupoid of isomorphisms between the fibers of the Lie group bundle $F$, otherwise 
known as the {\em frame groupoid} of $F$.

\begin{df} \label{def:xmlgpd}
A {\em crossed module} of Lie groupoids is a quadruple $xm = (F,\da,\Omega,\rho)$, where
$\Omega \gpd M$ is a Lie groupoid over $M,~F$ is a Lie group bundle on the same base,
$\da : F \rightarrow \Omega$ is a morphism of Lie groupoids over $M$, and where $\rho$ is
a representation of $\Omega$ on $F$, all such that
\begin{enumerate}
\item $\da(\rho(\xi,f)) = \xi \da(f) \xi^{-1}$ for all $(\xi,f) \in \Omega * F$;
\item $\rho(\da(f),f') = ff'f^{-1}$ for all $f,f' \in F$ with $\pi(f) = \pi(f')$;
\item $Im(\da)$ is a closed embedded submanifold of $\Omega$.
\end{enumerate}
\end{df}

The conditions of this definition show that $\im\da$ lies entirely in $I\Omega$ and is 
normal in $\Omega$. The normalcy of $\im\da$ then ensures that it is a Lie group bundle. The quotient 
$\Omega/\im\da$ therefore exists and is a Lie groupoid over $M$. This is called the {\em cokernel} of 
the crossed module and we usually denote it by $\Bar{\Omega}$. On the other hand, condition {\rm (ii)} 
ensures that $\ker\da$ lies in $ZF$. All this is described in figure 1.
\begin{figure}[!htb]
\begin{diagram}
\ker{\tau} & & & & & & \\
\dTeXembed & & \rho : \Omega * F & \rTeXto & F & \\
 F & & & & & \\
\dTeXonto<{\tau} & & & & & \\
Im(\tau) & \rTeXembed & \Omega & \rTeXonto^{\natural} & \Omega/Im(\tau) = \Bar{\Omega} & \\ 
\end{diagram}
\caption{}
\end{figure}

Notice that $\rho$ here also induces a representation of $\Bar{\Omega}$ on $\ker\da$, denoted by 
$\rho^{\ker\da}$, by
\[
\rho^{\ker\da}(\Bar{\xi},f) = \rho(\xi,f). 
\]
This is well defined because if we consider $\xi,\eta \in \Omega$ such that $\Bar{\xi} = \Bar{\eta}$ 
then there is a $f' \in F$ such that $\eta = \xi \cdot \da(f')$. So
\[
\rho(\eta,f) = \rho(\xi, \rho(\da(f'),f)) = = \rho(\xi,f'f(f')^{-1}) = \rho(\xi,f)
\]
since $\ker\da \subseteq ZF$.

Two crossed modules $(F,\da,\Omega,\rho)$ and $(F',\da',\Omega',\rho')$ are {\it equivalent} if there is a morphism of Lie groupoids $\theta : \Omega \to \Omega'$ such that $\theta \circ \da = \da, \sharp \circ \theta = \sharp'$ and $\rho' \circ \theta = \rho$. 
In the following the term crossed module will mean the relevant equivalence class and be denoted by $\langle F,\da,\Omega,\rho \rangle$. We regard a crossed module of Lie groupoids as a structure on the 
cokernel; if $\coker(\da) = \Bar{\Omega}$ we say that $\langle F,\da,\Omega,\rho \rangle$ is a 
{\em crossed module of $\Bar{\Omega}$ with $F$}. 

There are three special types of crossed modules worth noting on the groupoid level. First, crossed 
modules of Lie groupoids with trivial kernel. These are merely extensions of Lie groupoids 
$F \stackrel{\da}{\inj} \Omega \stackrel{\natural}{\surj} \Omega/F$ where $F \subseteq I\Omega$ is a normal 
subbundle of $I\Omega$ and the representation $\rho$ of $\Omega$ to $F$ is the restriction of the inner 
representation of $\Omega$ on $I\Omega$. namely, $\tau(\rho(\xi,f)) = I_{\xi}(f)$ for all 
$(\xi,f) \in \Omega * F$. The other two types are:

\begin{df}
A crossed module of Lie groupoids $\langle F,\da,\Omega,\rho \rangle$ over the manifold $M$ is called
\begin{enumerate}
\item a {\em coupling crossed module} if $\ker\da = ZK$;
\item a {\em pair crossed module} if $\coker\da = M \times M$.
\end{enumerate}
If both $\ker\da = ZK$ and $\coker\da = M \times M$, then $\langle F,\da,\Omega,\rho \rangle$ is called a 
{\em coupling pair} crossed module.
\end{df}

Any extension of Lie groupoids $F \stackrel{\iota}{\inj} \Omega \stackrel{\pi}{\surj} \Bar{\Omega}$ give rise to crossed modules of Lie groupoids. Namely the choice of a normal 
subbundle $N$ of $F$ which lies entirely in $ZF$ induces the crossed module 
$\langle F,\da,\Omega/N,\rho \rangle$ (of $\Bar{\Omega}$ with $F$), where $\da$ is the projection 
$F \rightarrow F/N$ and $\rho : \Omega/N * F \rightarrow F$ is the representation defined by 
\[
\rho(\langle \xi \rangle,f) = \xi \cdot \iota(f) \cdot \xi^{-1}
\]
The question whether every crossed module of Lie groupoids arises in this manner from an 
extension leads to the notion of an operator extension. 

\begin{df} \label{df:opextlalgd}
An {\em operator extension} of the crossed module of Lie groupoids $(F,\tau,\Omega,\rho)$
over the manifold $M$ with cokernel $\Bar{\Omega}$ is a pair $(\Hat{\Omega},\mu)$ where $\Hat{\Omega} \gpd M$ is  
a Lie groupoid extension of $\Bar{\Omega}$ by $F$ and $\mu : \Hat{\Omega} \rightarrow \Omega$ is a morphism
of Lie groupoids over $M$ such that:
\begin{enumerate}
\item The following diagram commutes:
$$
\begin{diagram}
          F          & \rTeXembed^{\iota} &  \Hat{\Omega}   & \rTeXonto & \Bar{\Omega} \\
\dTeXonto<{\tau} &                    & \dTeXonto<{\mu} &           &   \dequal    \\
       \im\da        &     \rTeXembed     &    \Omega       & \rTeXonto & \Bar{\Omega} \\
\end{diagram}
$$
\item $\iota(\rho(\mu(\Hat{\xi}),f)) = \Hat{\xi} \cdot \iota(f) \cdot \Hat{\xi}^{-1}$ for all
$(\Hat{\xi},f) \in \Hat{\Omega} * F$.
\end{enumerate}
\end{df}
\begin{df}
Let $(F,\tau,\Omega,\rho)$ be a pair crossed module of Lie groupoids over the manifold
$M$. Two operator extensions $(\Hat{\Omega},\mu)$ and $(\Hat{\Omega}',{\mu}')$ are called
{\em equivalent} if there is an isomorphism of Lie groupoids $\kappa : \Hat{\Omega} \rightarrow \Hat{\Omega}'$ such that ${\mu}' \circ \kappa = \mu$.
\end{df} 
 
The obstruction associated with a pair crossed module of Lie groupoids was given by  Mackenzie in \cite{Mackenzie:1989}. We will show in Section 4 that in order to understand the obstruction associated with a general crossed module of Lie groupoids it suffices to understand the 
obstruction in the particular case of pair crossed modules. Let us recall in brief the construction of the obstruction from \cite{Mackenzie:1989}. 

A pair crossed module of Lie groupoids $\langle F,\da,\Omega,\rho \rangle$ is, as we discussed earlier, a crossed module of $M \times M$ with $F$. This means that the gauge group bundle of the Lie groupoid $\Omega$ is the image of $\da$. The groupoid $\Omega$ can therefore be written as an extension of groupoids in the form 
\[
\im\da \inj \Omega \stackrel{(\beta,\alpha)}{\surj} M \times M.
\]
Fix an element $x_{0} \in M$ and denote the Lie group $F_{x_{0}}$ by $H$. Choose an open simple cover $\{ U_{i} \}_{i \in I}$ of $M$. We write $U_{ij}$ for the intersection of two open sets $U_{i}$ and $U_{j}$, also $U_{ijk}$ for the intersection of three open sets, etc.  
Let $\{ s_{ij} : U_{ij} \rightarrow \da(H) \}_{i,j \in I}$ be a cocycle of transition functions for the Lie 
groupoid $\Omega \gpd M$ and $\{ \Hat{s}_{ij} : U_{ij} \rightarrow H \}_{i,j \in I}$ be smooth lifts of the transition functions to $H$, such that $s_{ij} = \da \circ \Hat{s}_{ij}$. Now consider the failure of these lifts to form a cocycle
\[
e_{ijk} : U_{ijk} \rightarrow H,~e_{ijk} = \Hat{s}_{jk} \cdot \Hat{s}_{ik}^{-1} \cdot \Hat{s}_{ij}.
\]
It follows from the fact that the $s_{ij}$'s form a cocycle that this function takes values in $ZH$. Therefore 
it defines a class $[e] \in \check{H}^{2}(M,ZH)$. This class depends neither on the choice of cocycle for 
$\Omega$ nor from the choice of lifts for this cocycle. Of course it is zero if and only if the lifts 
$\Hat{s}_{ij}$ form a cocycle and in this case they define a Lie groupoid $\Hat{\Omega} \gpd M$. It is proven 
in \cite{Mackenzie:1989} that $\Hat{\Omega}$ is an operator extension for the crossed module. 
This element is called the {\em obstruction} of the crossed module and we denote it by 
$\Obs\langle F,\da,\Omega,\rho \rangle$. In \cite{Mackenzie:1989} it was also shown that if 
$\Obs\langle F,\da,\Omega,\rho \rangle = 0$ then the equivalence classes of operator extensions of the crossed 
module $\langle F,\da,\Omega,\rho \rangle$ are classified by $\check{H}^{1}(M,ZH)$. 

\section{Crossed modules of Lie algebroids}

Let us discuss the differentiation of a crossed module as above. 
Given a crossed module of Lie groupoids $\langle F,\da,\Omega,\rho \rangle$ over $M$ it is well known that the 
Lie group bundle $F$ differentiates to a Lie algebra bundle $F_{*}$ over $M$, the Lie groupoid $\Omega$ to 
a Lie algebroid $A\Omega$ over $M$ and the morphism $\da$ to a morphism of Lie algebra bundles 
$\da_{*} : F_{*} \rightarrow L\Omega \subseteq A\Omega$. The part that needs some attention is the 
differentiation of the representation $\rho$. For every $\xi \in \Omega$ the map 
$\rho(\xi) : F_{\alpha(\xi)} \rightarrow F_{\beta(\xi)}$ is a Lie group isomorphism. The Lie functor then 
shows that it differentiates to an isomorphism of Lie algebras 
$(\rho(\xi))_{*} : (F_{\alpha(\xi)})_{*} \rightarrow (F_{\beta(\xi)})_{*}$. Denoting $\Phi(F_{*}) \gpd M$ the 
Lie groupoid of isomorphisms between the fibers of the Lie algebra bundle $F_{*}$ (otherwise known as the frame 
groupoid of $F_{*}$), we get a well defined morphism of Lie groupoids
\[
\Tilde{\rho} : \Omega \rightarrow \Phi(F_{*}), \xi \mapsto (\rho(\xi))_{*}
\] 
Now apply the Lie functor to $\Tilde{\rho}$ to get the morphism of Lie algebroids 
$\rho_{*} : A\Omega \rightarrow CDO[F_{*}]$. This is the representation $\rho$ differentiates to. 
\begin{lem}\label{lemma:diffxm}
\begin{enumerate}
\item $\rho_{*}(\da_{*}(V)) = \ad_{V}$ for all $V \in F_{*}$ and
\item $\da_{*} \circ \rho_{*}(X) = \ad_{X} \circ \da_{*}$ for all $X \in A\Omega$.
\end{enumerate}
\end{lem}
  
\pf~Using the definitions of $\rho_{*}$ and $\tau_{*}$ we have:
\[
\rho_{*}(\tau_{*}(V)) = \rho_{*}(T_{e_{x}}\tau(V)) = (T_{1_{x}}\Tilde{\rho} \circ
T_{e_{x}}\tau)(V) = T_{e_{x}}(\Tilde{\rho} \circ \tau)(V)
\]
for any $V \in (F_{*})_{x}$ and $x \in M$. On the other hand, for all $f \in F$ we have
\[
(\Tilde{\rho} \circ \tau)(f) = T_{e_{\pi(f)}}(\rho(\tau(f))) = T_{e_{\pi(f)}}(I_{f}) =
Ad_{f}
\]
So,
\[
\rho_{*}(\tau_{*}(V)) = T_{e_{x}}(\Tilde{\rho} \circ \tau)(V) = T_{e_{x}}Ad(V) = ad_{V}.
\]
For {\rm (ii)} we know that $\tau \circ \rho(\xi) = I_{\xi} \circ \tau$ for
all $\xi \in \Omega$. Therefore,
\[
T_{e_{\alpha(\xi)}}(\tau \circ \rho(\xi)) = T_{e_{\alpha(\xi)}}(I_{\xi} \circ \tau) =
\Ad_{\xi} \circ T_{e_{\alpha(\xi)}}\tau 
\Rightarrow T_{e_{\beta{\xi}}}\tau \circ \Tilde{\rho}(\xi) = \Ad_{\xi} \circ
T_{e_{\alpha(\xi)}}\tau.
\]
By differentiating the last equality and using the fact that $T_{e_{x}}\da$ is linear,
therefore it is its own derivative, we get $\da \circ \rho_{*} (X) = \ad_{X} \circ \da_{*}$. \boom

This naturally leads to the following definition.
\begin{df}
\label{df:xm}
A {\em crossed module of Lie algebroids} over the manifold $M$ is 
a quadruple $(K,\da,A,\rho)$ where $K \to M$ is a Lie algebra bundle, 
$A\to M$ is a transitive Lie algebroid, $\da \co K \to A$ is a morphism 
of Lie algebroids and $\rho \co A \to \CDO[K]$ is a representation of 
$A$ in $K$ such that:
\begin{enumerate}
\item $\rho(\da(V))(W) = [V,W]$ for all $V,W \in \Gamma K$ and
\item $\da(\rho(X)(V)) = [X,\da(V)]$ for all 
$X \in \Gamma A,\ V \in \Gamma K$.
\end{enumerate}
\end{df}

Since $\da$ is a morphism of Lie algebroids, we have 
$a \circ \da = 0$, where $a$ is the anchor of $A$. So
$\im(\da)$ lies entirely in the adjoint bundle $L$ of $A$.
Regarding $\da$ temporarily as a morphism of Lie algebra bundles, 
condition (ii) is equivariance with respect to $\rho$ and the
adjoint action of $A$ on $L$. It is therefore of locally
constant rank (see \cite{KCHM:new}, discussion after I, 3.3.13 and 6.5.11), and so 
has a kernel Lie algebra bundle which we denote $\ker{\da}$. 
Condition (i) now ensures that $\ker{\da}$ lies in $ZK$. 
Likewise, the quotient Lie algebroid $A/\im(\da)$ exists and is 
a Lie algebroid over $M$ (see \cite[I\S4.4]{KCHM:new}). 
This is called the {\em cokernel} of the crossed module, and
we usually denote it $\Bar{A}$.  
All this is described in figure 2:
\begin{figure}[htp]
\begin{diagram}
\ker{\da} & & & & & & \\
\dTeXembed & & \rho : A & \rTeXto & \CDO[K] & \\
 K & & & & & \\
\dTeXonto<{\da} & & & & & \\
\im(\da) & \rTeXembed & A & \rTeXonto^{\natural} & A/\im(\da) = \Bar{A} & \\
\end{diagram}
\caption{}
\end{figure}

Notice that $\rho$ induces a representation of $\Bar{A}$ on the
vector bundle $\ker{\da}$, denoted $\rho^{\ker{\da}}$, by
$$
\rho^{\ker{\da}}(\Bar{X})(V) = \rho(X)(V).
$$
This is well defined because if we consider
$X,Y \in \Gamma A$ such that $\Bar{X}  = \Bar{Y} \in \Bar{A}$, 
then there is a $W \in \Gamma K$ such that $X = Y + \tau(W)$. So
$$
\rho(X)(V) = \rho(Y)(V) + \rho(\tau(W))(V) = 
= \rho(Y)(V) + [W,V] = \rho(Y)(V)
$$
since $\ker{\da}\subseteq ZK$. 

Throughout the rest of the paper we will be working with crossed modules 
with fixed cokernel, as well as fixed $K$ and $\ker{\tau}$. Two such 
crossed modules $(K,\tau,A,\rho)$ and $(K,\tau',A',\rho')$ 
are {\em equivalent} if there is a morphism of Lie algebroids
$\theta \co A \to A'$ such that $\theta \circ \tau = \tau$,
$\natural \circ \theta = \natural'$ and $\rho' \circ \theta = \rho$. 
The 3--lemma then shows that every such morphism is an
isomorphism of Lie algebroids. From now on, every time we mention 
a crossed module of Lie algebroids we will refer to its equivalence 
class and denote it by $\langle K,\tau,A,\rho \rangle$.

Again, there are three special types of Lie algebroid crossed modules worth noting. 
The ones with with trivial kernel are merely extensions of Lie algebroids
$K \stackrel{\tau}{\inj} A \stackrel{\natural}{\surj} A/K$ where $K\subseteq L$
is an ideal of $A$, and the representation $\rho$ of
$A$ in $K$ is the restriction to $K$ of the adjoint representation of $A$ on $L$. 
Namely, $\tau(\rho(X)(V)) = \ad_{X}(\tau(V))$ for all $X \in \Gamma A$ and $V \in \Gamma K$. We will also need the following ones:

\begin{df}
A crossed module of Lie algebroids $\langle K,\tau,A,\rho \rangle$ 
over the manifold $M$ is called
\begin{enumerate}
\item a {\em coupling crossed module} 
if $\ker{\tau} = ZK$;
\item a {\em pair crossed module} if $\coker\da = TM$.
\end{enumerate}
If both $\ker{\tau} = ZK$ and $\coker\da = TM$, then 
$\langle K,\tau,A,\rho \rangle$ is called a 
\emph{coupling pair crossed module.}
\end{df}

We usually regard a coupling crossed module as a structure on the cokernel;
if $\coker(\da) = \Bar{A}$ we say that $\langle K,\tau,A,\rho \rangle$ is
a \emph{coupling crossed module of $\Bar{A}$ with $K$}.

\subsection{Some examples}

Let us discuss some examples of Lie algeberoid crossed modules.

\begin{enumerate}
\item Take any principal bundle $P(M,G,\pi)$ and consider its Atiyah sequence
\[
\frac{P \times \gog}{G} \stackrel{j}{\inj} \frac{TP}{G} \stackrel{\pi_{*}}{\surj} TM.
\]
(the action of $G$ on $\gog$ is the adjoint). This is naturally a Lie algebroid. Recall the identification $j(\frac{P \times \gog}{G}) = \frac{T^{\pi}P}{G}$ induced by the map $\Tilde{j} : P \times \gog \to TP$, $(u,X) \mapsto T_{(u,1)}m(0_{u},X_{1})$, where $m : P \times G \to P$ is $(u,g) \mapsto ug$.

Now quotient the Lie algebroid $\frac{TP}{G}$ by the center of the Lie algebra bundle $\frac{P \times \gog}{G}$. This quotient may be identified with the Lie subalgebroid $\ad (\frac{TP}{G})$ of 
$\CDO [ \frac{P \times \gog}{G} ]$, namely the image of the representation $\ad : \frac{TP}{G} \to \CDO [ \frac{P \times \gog}{G} ]$ given by
\[
\ad_{\langle X \rangle}(\langle (u,V) \rangle) = j^{-1}([ X,\Tilde{j}^{-1}(u,V) ])
\]
Then $\langle \frac{P \times \gog}{G}, \ad\mid_{\frac{P \times \gog}{g}}, \ad(\frac{TP}{G}),\ad \rangle$ is the induced crossed module.

\item Let $(M,\omega)$ be a symplectic manifold. Then the vector bundle $A = TM \oplus (M \times \R)$ becomes a Lie algebroid with bracket
\[
[ X \oplus V,Y \oplus W ] = [ X,Y ] \oplus \{ X(W) - Y(V) - \omega(X,Y) \}
\]
and anchor the first projection. This algebroid naturally induces the crossed module $\langle M \times \R, \da, A, \rho \rangle$ where $\da : M \times \R \to A$ is the natural inclusion $V \mapsto 0 \otimes V$ and $\rho : A \to \CDO [ M \times \R ]$ is the representation $\rho(X \oplus V)(W) = 0 \oplus X(W)$

\item A non-abelian generalisation of a symplectic form was given on \cite[8.3.9]{KCHM:new}, where the classical Weil lemma was extended. Namely, let $L$ be a Lie algebra bundle over $M$ together with a connection $\nabla$ such that $\nabla_{X}[ V,W ] = [ \nabla_{X}V,W ] + [ V,\nabla_{X}W ]$. Then any $L$-valued 2-form $R$ on $M$ such that $R_{\nabla} = \ad \circ R$ (where $R_{\nabla}$ stands for the curvature of $\nabla$) and $\nabla(R) = 0$ endows the sections of the vector bundle $A = TM \oplus L$ with the Lie bracket
\[
[ X \oplus V,Y \oplus W ] = [ X,Y ] \oplus \{ \nabla_{X}W - \nabla_{Y}V + [V,W] - R(X,Y) \}
\]
which makes $A$ a Lie algebroid. Now consider the quotient bundle $A/ZL$. This may be identified with the Lie subalgebroid $\ad(A)$ of $\CDO [ L ]$, namely the image of the representation $\ad : A \to \CDO [ L ]$ defined by $\ad_{X \oplus V} 0 \oplus W = [ X \oplus V, 0 \oplus W ]$. The identification is $(X \oplus V) + ZL \mapsto \ad_{X \oplus V}$. It is straightforward that the induced crossed module is $\langle L \ad \mid_{L}, \ad(A), \ad \rangle$.

Notice that the algebroid $\ad(A)$ is integrable as a Lie subalgebroid of $\CDO [ L ]$.
\end{enumerate}

\subsection{Equivalence with couplings}

The notion of coupling crossed module is equivalent to the concept of
coupling of Lie algebroids, introduced in \cite[I, \S 7.2]{KCHM:new}
as the Lie algebroid form of the notion of ``abstract kernel'' in the sense 
of MacLane \cite{MacLane:1975}. Let us discuss how this equivalence is established.

Consider a Lie algebra bundle $K$ over the manifold $M$. The adjoint bundle of
$\CDO[K]$ is $\Der(K)$, the derivations of $K$, and 
$\ad(K) = \im(\ad \co K \to \Der(K))$ is a Lie subalgebra bundle of $\Der(K)$, 
and an ideal of $\CDO[K]$. We denote the quotient Lie algebroid $\CDO[K]/\ad(K)$ 
by $\OutDO[K]$, and call elements of
$\Gamma \OutDO[K]$ {\em outer covariant differential operators} on $K$.

\begin{df}
\label{df:coupling}
A {\em coupling} of the Lie algebroid $\Bar{A}$ with the Lie algebra bundle 
$K$ (both over the same manifold $M$) is a morphism of Lie algebroids 
$\Xi \co \Bar{A} \to \OutDO[K]$.
\end{df}

Fix a coupling $\Xi$ of the Lie algebroid $\Bar{A}$ with the Lie algebra bundle 
$K$. Since the map $\natural \co \CDO[K] \to \OutDO[K]$ is a surjective 
submersion as a map of vector bundles over $M$, there is a vector bundle morphism
$\nabla \co \Bar{A} \to \CDO[K],\ X \mapsto \nabla_{X}$, such that 
$\natural \circ \nabla = \Xi$. We call $\nabla$ a {\em Lie derivation law 
covering} $\Xi$.

Let $\nabla$ be any such Lie derivation law. Then for $X \in \Gamma \Bar{A}$ 
the operator $\nabla_{X} \co \Gamma K \to \Gamma K$ restricts to 
$\Gamma ZK \to \Gamma ZK$, for if $Z \in \Gamma ZK$ and $V \in \Gamma K$ then
$$
[V,\nabla_{X}(Z)] = \nabla_{X}([V,Z]) - [\nabla_{X}(V),Z] = \nabla_{X}(0) - 0 = 0,
$$
since $Z$ is central. Further, the restriction is independent of the choice
of $\nabla$; write $\rho^{\Xi}$ for the restriction of 
$\nabla_{X}$ to $\Gamma ZK \to \Gamma ZK$. Then $\rho^{\Xi}$ defines 
a vector bundle map $\Bar{A} \to \CDO(ZK)$ which is easily seen to be a 
Lie algebroid morphism; that is, $\rho^\Xi$ is a representation of $\Bar{A}$
on $ZK$, called the {\em central representation} of $\Xi$.

Now we can proceed to prove the equivalence of couplings in the sense 
of \ref{df:coupling} with coupling crossed modules of Lie algebroids. 

Consider first a coupling crossed module $\langle K,\tau,A,\rho \rangle$ 
of the Lie algebroid $\Bar{A}$ with the Lie algebra bundle $K$. Condition (i)
of \ref{df:xm} shows that $\rho$ sends $\im(\da)$ to $\ad(K)\subseteq\CDO[K]$
and $\rho$ therefore descends to a morphism $\Xi^{\rho}\co \Bar{A} \to \OutDO[K]$ 
as in the diagram
$$
\begin{diagram}
A & \rto^{\rho} & \CDO[K] \\
\dto<{} &  & \dto>{} \\
\Bar{A} & \rto_{\Xi^{\rho}} & \OutDO[K]
\end{diagram}
$$
where the two vertical maps are the natural projections. This $\Xi^\rho$ is 
the coupling corresponding to $\langle K,\tau,A,\rho \rangle$. Note that
equivalent coupling crossed modules induce the same coupling. It is also easy 
to see that the representation of $\Bar{A}$ on $ZK$ induced by $\Xi^\rho$ is 
equal to the representation induced directly from $\rho$ as in the passage 
following \ref{df:xm}. 

For the construction of the coupling crossed module corresponding to a coupling
we use the construction principle of \cite[II~7.3.7]{KCHM:new}. Take a 
coupling $\Xi \co \Bar{A} \to \OutDO[K]$ of the Lie algebroid $\Bar{A}$ with 
the Lie algebra bundle $K$. Choose a Lie derivation law 
$\nabla \co \Bar{A} \to \CDO[K]$ covering $\Xi$. This is an
anchor--preserving vector bundle morphism, and so its curvature is a 
well defined map
$$
R_{\nabla} \co \Bar{A} \oplus \Bar{A} \to \CDO[K],\qquad
X\oplus Y\to \nabla_{[X, Y]} - [\nabla_X, \nabla_Y].
$$ 
Since $\natural \circ \nabla = \Xi$ is a morphism of Lie algebroids, 
it follows that $\natural \circ R_{\nabla} = 0$ and so $R_{\nabla}$
takes values in $\ad(K) \subseteq \Der(K)$.

Define a map $\Bar{\nabla} \co \Bar{A} \to \CDO[\ad(K)]$ by
$$
\Bar{\nabla}_{X}(\ad_{V}) = \ad_{\nabla_{X}(V)}
$$
for all $X \in \Gamma \Bar{A}$ and $V \in \Gamma K$. This is also an anchor 
preserving morphism, so its curvature is a well defined map 
$R_{\Bar{\nabla}} \co \Bar{A} \oplus \Bar{A} \to \CDO[\ad(K)]$. It is easily 
verified that $R_{\Bar{\nabla}} = \ad \circ R_{\nabla}$ and
$\Bar{\nabla}(R_{\nabla}) = 0$. Moreover, the map 
$\natural \circ \Bar{\nabla} \co \Bar{A} \to \OutDO[\ad(K)]$ has zero curvature. 

Now \cite[II~7.3.7]{KCHM:new} shows that the formula
$$
[X \oplus \ad_{V},Y \oplus \ad_{W}] = 
[X,Y] \oplus \{ \ad_{\nabla_{X}(W)} - \ad_{\nabla_{Y}(V)} +
        \ad_{[V,W]} - R_{\nabla}(X,Y) \}
$$
defines a Lie bracket on $\Gamma(\Bar{A} \oplus \ad(K))$ which makes 
$\Bar{A} \oplus \ad(K)$ a Lie algebroid
over $M$. Denote $\Bar{A} \oplus \ad(K)$ by $A$. 
Define $\tau \co K \to A$ and $\rho \co A \to \CDO[K]$ by
$$
\tau(V) = 0 \oplus \ad_{V}, 
\qquad\qquad
\rho(X \oplus \ad_{V}) = \nabla_{X}(W) + [V,W]
$$
for all $X \in \Gamma \Bar{A}$ and $V,W \in \Gamma K$. 
These are both morphisms of Lie algebroids and the remaining steps in the
following proof are straightforward. 

\begin{prop}
The Lie algebroid $A$ just defined, together with $\da$ and $\rho$, constitute
a coupling crossed module for $\Bar{A}$, which induces the given $\Xi$. 
\end{prop}

\subsection{Lifting crossed modules of Lie algebroids}

A transitive Lie algebroid $A$ over $M$ is in fact an extension of the tangent bundle $TM$ 
by its adjoint bundle $L$. From this point of view, an arbitrary exptension of Lie algebroids 
$K \stackrel{\iota}{\inj} A \stackrel{\pi}{\surj} \Bar{A}$ over a manifold $M$ gives rise 
to a crossed module of Lie algebroids once we choose an ideal $I$ of $ZK$, the center of the 
Lie algebra bundle $K$. The construction is described in figure 3:
\begin{figure}[!htb]
$$
\begin{diagram}
I & \requal & I &  & \\
\dTeXembed & & \dTeXembed &  & \\
K & \rTeXembed^{\iota} & A & \rTeXonto^{\pi} & \Bar{A} \\
\dTeXonto<{\da^{I}} &  & \dTeXonto<{\natural} &  & \dequal \\
\im(\da^{I}) & \rTeXembed^{\iota^{\natural}} & A/\iota(I) & \rTeXonto^{\pi^{\natural}} & \Bar{A} 
\end{diagram}
$$  
\caption{}
\end{figure}

Here $\da^{I}$ is the quotient map $K \rightarrow A/\iota(I)$ and the other maps are 
$\iota^{\natural}(\da^{I}(V)) = \langle \iota(V) \rangle$ for all $V \in K$ and 
$\pi^{\natural}(\langle X \rangle) = \pi(X)$ for all $X \in A$. These maps are well defined 
because the sequence $K \stackrel{\iota}{\inj} A \stackrel{\pi}{\surj} \Bar{A}$ is exact. Now the 
representation $\rho^{I} : A/\iota(I) \rightarrow CDO[K]$ defined by 
\[
\iota(\rho^{I}(\langle X \rangle)(V)) = [X,\iota(V)]
\]
is well defined because we assumed $I$ to be an ideal of $ZK$ and it
makes the quadruple $\langle K,\da^{I},A/\iota(I),\rho \rangle$ a crossed module of Lie algebroids. 
The question whether every crossed module arises from an extension of Lie algebroids gives rise to the notion of an operator extension. 
\begin{df}\label{df:opextlalgds}
Let $\xm = \langle K,\tau,A,\rho \rangle$ be a crossed module of Lie algebroids of$\Bar{A}$ with $K$. 
An {\em operator extension} of $\xm$ is a pair
$(K \stackrel{\iota'}{\inj} \hat{A} \stackrel{\pi}{\surj} \bar{A},\mu_{*})$ of an extension of Lie
algebroids together with a morphism of Lie algebroids $\mu : \hat{A} \rightarrow A$ which is a
surjective submersion such that:
\begin{enumerate}
\item The following diagram commutes:
$$
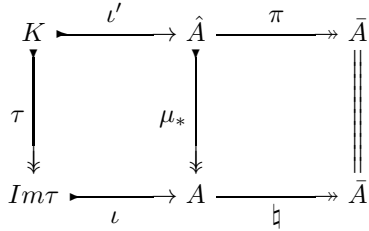
\begin{diagram}
            K          & \rTeXembed^{\iota'} &      \hat{A}      &    \rTeXonto^{\pi}   & \bar{A} \\
\dTeXembedd<{\tau} &                     & \dTeXembedd<{\mu_{*}} &                      & \dequal \\
       Im\tau      &  \rTeXembed_{\iota} &          A        & \rTeXonto_{\natural} & \bar{A} \\
\end{diagram}
$$
\item For all $\hat{X}~\in~\Gamma A,~V~\in~\Gamma K$ we have:
\[
\iota(\rho(\mu_{*}(\hat{X}))(V)) = [\hat{X},\iota(V)]
\]
\end{enumerate}
\end{df} 

\begin{df}
The operator extensions $(K\stackrel{\iota'_{1}}{\inj} \hat{A}_{1} \stackrel{\pi_{1}}{\surj}
\bar{A},\mu_{*}^{1})$ and $(K\stackrel{\iota'_{2}}{\inj} \hat{A}_{2} \stackrel{\pi_{2}}{\surj}
\bar{A},\mu_{*}^{2})$ of the crossed module of Lie algebroids $\xm = \langle K,\tau, A,\rho
\rangle$ are {\em equivalent} if there is a Lie algebroid morphism
$\kappa_{*} : \hat{A}_{1} \rightarrow \hat{A}_{2}$ such that $\mu_{*}^{2} \circ \psi = \mu_{*}^{1}$ and the
following diagram commutes:
$$
\begin{diagram}
   K    & \rTeXembed^{\iota'_{1}} &   \hat{A}_{1}  & \rTeXonto^{\pi_{1}} & \bar{A} \\
\dequal &                         & \dTeXto<{\kappa_{*}} &                     & \dequal \\
   K    & \rTeXembed^{\iota'_{2}} &   \hat{A}_{2}  & \rTeXonto^{\pi_{2}} & \bar{A} \\
\end{diagram}
$$
We denote $\Opext(\xm)$ the set of equivalence classes of operator extensions of the crossed module
$\xm$.
\end{df}

\begin{rmk}[The lifting obstruction]
The obstruction to the existence of an operator extension for a crossed module of Lie algebroids is the same as the one given in \cite[II\S7.3]{KCHM:new} for couplings. This obstruction is a certain class in ${\cal{H}}^{3}(TM,\rho^{\da},ZK)$ (Lie algebroid cohomology). When it vanishes the equivalence classes of operator extensions (lifts) are classified by ${\cal{H}}^{2}(TM,\rho^{\da},ZK)$.
\end{rmk}

\section{The integrability of transitive Lie algebroids via crossed modules}

This section provides an alternative way to obtain the integrability obstruction of a transitive Lie algebroid $A \to M$. We prove that the integrability obstruction given by Mackenzie in \cite{KCHM:new} coincides with the lifting obstruction of a certain crossed module of Lie groupoids naturally associated with $A$. Let us begin with the following result which clarifies the nature of the lifting obstruction. It actually motivates our approach, since it shows that the lifting obstruction is a cohomology class of the same type as the integrability obstruction given by Mackenzie.

\begin{prop}\label{constant:opextobstr}
If a pair crossed module of Lie groupoids differentiates to a coupling crossed module 
of Lie algebroids for which the obstruction class vanishes, then the obstruction class 
for the crossed module of Lie groupoids takes values in \v{C}ech cohomology with coefficients in constant functions. 
\end{prop}

\pf Consider a pair crossed module of Lie groupoids $\langle F,\da,\Omega, \rho \rangle$. If the obstruction 
class of the coupling $\langle F_{*},\da_{*},A\Omega,\rho_{*} \rangle$ vanishes then it has an operator extension, 
i.e. there exists a (transitive) Lie algebroid $\Hat{A}$ over $M$ and a morphism of Lie algebroids 
$\mu_{*} : \Hat{A} \rightarrow A\Omega$ which is a surjective submersion such that the diagram in figure 4 commutes:
\begin{figure}[!hbt]
$$
\begin{diagram}
F_{*} & \rTeXembed^{\iota'} & \Hat{A} & \rTeXonto^{\pi} & TM \\
\dTeXonto<{\da_{*}} &  & \dTeXonto<{\mu_{*}} & & \dequal \\
\im \da_{*} & \rTeXembed_{\iota} & A\Omega & \rTeXonto_{\natural} & TM. 
\end{diagram}
$$
\caption{}\label{operatorextension}
\end{figure}

Let $H$ denote the fiber of the Lie group bundle $F \rightarrow M$ and $\hoh$ its Lie algebra. Then 
the vertex groups of the Lie groupoid $\Omega \gpd M$ are isomorphic to $\da(H)$. Take a simple open cover 
$\{ U_{i} \}_{i \in I}$ of $M$. Choose a section atlas $\{ s_{ij} : U_{ij} \rightarrow \da(H) \}_{i,j \in I}$ for 
the Lie groupoid $\Omega \gpd M$ over this cover. Consider the family of Maurer-Cartan forms 
$\chi_{ij} : TU_{ij} \rightarrow U_{ij} \times\da_{*}(\hoh)$ defined as $\chi_{ij} = \Delta(s_{ij})$ and the 
cocycle $\alpha_{ij} = \Ad_{s_{ij}}$ with values in $\Aut(H)$. Here $\Delta$ denotes the Darboux derivative, 
otherwise known as the right-derivative of functions with values in a Lie group. These define a system of transition 
data $(\chi,\alpha)$ in the sense of \cite[II 8.2.5]{KCHM:new}. This means that this pair is compatible 
in the sense 
\begin{eqnarray}\label{compatibility}
\Delta(\alpha_{ij}) = \ad \circ \chi_{ij}
\end{eqnarray} 
and it satisfies the cocycle condition 
\begin{eqnarray}\label{cocycle}
\chi_{ik} = \chi_{ij} + \alpha_{ij}(\chi_{jk}) 
\end{eqnarray}
for all $i,j,k \in I$. In fact, it is proven in \cite[II 8.2.5]{KCHM:new} that systems of transition 
data classify transitive Lie algebroids. This and the commutativity of the diagram in figure \ref{operatorextension} show 
that there exists a system of transition data $(\Hat{\chi},\Hat{\alpha})$ with values in $\hoh$ for the Lie 
algebroid $\Hat{A}$ such that $\da_{*} \circ \Hat{\chi_{ij}} = \chi_{ij}$ and 
$\da_{*} \circ \Hat{\alpha}_{ij} = \alpha_{ij}$. The $\Hat{\chi}_{ij}$'s are Maurer-Cartan forms, so they 
integrate uniquely to smooth functions $\Hat{s}_{ij} : U_{ij} \rightarrow H$ such that 
$\Hat{\chi}_{ij} = \Delta(\Hat{s}_{ij})$. Following the same steps as 
in \cite[II \S8.3]{KCHM:new} it is proven that the compatibility condition \ref{compatibility} for the 
system of transition data $(\Hat{\chi},\Hat{\alpha})$ gives 
\[
\Hat{\alpha}_{ij} = \Ad_{\Hat{s}_{ij}}.
\] 
The system $(\Hat{\chi},\Hat{\alpha})$ satisfies the cocycle condition \ref{cocycle}. This gives
\[
\Delta(\Hat{s}_{ik}) = \Delta(\Hat{s}_{ij}) + \Ad_{\Hat{s}_{ij}}(\Delta(\Hat{s}_{jk})) 
\Rightarrow \Delta(\Hat{s}_{ik}) = \Delta(\Hat{s}_{ij} \cdot \Hat{s}_{jk}).
\]
A uniqueness argument now shows that there is a constant $c_{ijk} \in H$ such that 
$e_{ijk} = \Hat{s}_{jk} \cdot \Hat{s}_{ik}^{-1} \cdot \Hat{s}_{ij} = c_{ijk}$. In fact, $c_{ijk}$ lies in the 
center of the Lie group $H$ because $e_{ijk}$ takes values exactly there. So, the obstruction class of the 
pair crossed module of Lie groupoids $\langle F,\da,\Omega,\rho \rangle$ takes values in \v{C}ech cohomology with constant coefficients.~\boom

\begin{rmk}
A crossed module of Lie algebroids $\langle K,\da,A,\rho \rangle$ integrates to a crossed module of Lie 
groupoids if the Lie algebroid $A$ is integrable. A proof of this is given in Section 4.
Even when a crossed module integrates though, it does not follow that its operator extensions (if there are any) also integrate. Examples of such a situation are 
the non-integrable transitive Lie algebroids. Every such algebroid $A$ over a manifold $M$ is a Lie algebroid 
extension $L \inj A \surj TM$. The crossed module induced by the choice of an ideal $I \subseteq ZL$ may be integrable, 
but the crossed module of Lie groupoids it would integrate to can not have an operator extension. If an operator 
extension on the groupoid level existed, then it would have to differentiate to $L \inj A \surj TM$ and in this case 
the Lie algebroid $A$ would be integrable, which is a contradiction. More particularly, one may consider the example of the Lie algebroid associated to a symplectic manifold $(M,\omega)$ where the 2-form $\omega$ is not integral.
\end{rmk}

We may now continue with the reformulation of Mackenzie's integrability obstruction. Consider a transitive Lie algebroid $L \stackrel{\iota}{\inj} A\stackrel{\alpha_{A}}{\surj} TM$. The ideal $ZL$ induces the coupling crossed module $\langle L,\natural,A/ZL,\rho_{*} \rangle$, where $A/ZL$ is 
a Lie algebroid over $M$ with Lie bracket 
\[
[X + ZL,Y + ZL]_{A/ZL} = [X,Y]_{A} + ZL  
\] 
for all $X,Y \in \Gamma A$. Its anchor map is $q_{A/ZL}(X + ZL) = q_{A}(X)$, therefore it can be written as an extension 
of Lie algebroids in the form 
\[
L/ZL \inj A/ZL \stackrel{q_{A/ZL}}{\surj} TM. 
\]
The map $\natural : L \rightarrow A/ZL$ is of course the quotient projection and its image is the adjoint bundle$L/ZL$. Notice 
that the Lie algebroid $A/ZL$ can be identified canonically with the Lie subalgebroid $\ad(A)$ of $\CDO[L]$. This is the image 
of the adjoint representation $\ad : A \rightarrow CDO[L]$. Namely, every 
$X +ZL \in \Gamma A/ZL$ defines the operator $\ad_{X} : \Gamma L \rightarrow \Gamma L$ by $\ad_{X}(V) = [X,\iota(V)]_{A}$ 
for every $V \in \Gamma L$. This element is well defined because if $Y \in \Gamma A$ is another representative of the class 
$X + ZL$ then there exists an element $W \in \Gamma ZL$ such that $Y = X +\iota(W)$. Then
\[
\ad_{Y}(V) = [Y,\iota(V)]_{A} = [X,\iota(V)]_{A} + [\iota(W),\iota(V)]_{A} = [X,\iota(V)]_{A} = \ad_{X}(V),
\]
since $W \in \Gamma ZL$. On the other hand, every element $\ad_{X} \in \ad(A)$ can be canonically identified with $X + ZL \in A/ZL$. 
A similar argument shows that $L/ZL$ can be identified canonically with $\ad(L)$, the bundle of inner automorphisms of the 
fibers of $L$.  

Now that we have established this identification, it is easier to regard the representation 
$\rho_{*} : A/ZL \rightarrow \CDO[L]$ as the natural inclusion of algebroids $\ad(A) \inj \CDO[L]$ and the quotient map $\natural$ 
as $\ad : L \rightarrow \ad(L)$. All this is described in figure 5:
\begin{figure}
\begin{diagram}
ZL & \requal & ZL &  &  &\rho_{*} : \ad(A) \rightarrow \CDO[L]\\
\dTeXembed & & \dTeXembed &  &  &\\
L & \rTeXembed^{\iota} & A & \rTeXonto^{\pi} & TM &\\
\dTeXonto<{\ad} &  & \dTeXonto<{\ad} &  & \dequal &\\
\ad(L) & \rTeXembed^{\iota^{\natural}} & \ad(A) & \rTeXonto^{\pi^{\natural}} & TM & \\
\end{diagram}
\caption{}
\end{figure}

The coupling induced by $A$ can now be written as $\langle L,\ad,\ad(A),\rho_{*} \rangle$. 
The Lie algebroid $\CDO[L]$ integrates to the frame groupoid $\Phi[L] \gpd M$ (see \cite[I 3.6.6]{KCHM:new}). Therefore, 
$\ad(A)$ also integrates as a Lie subalgebroid of $\CDO[L]$. The Lie groupoid it integrates to is denoted by $\Int(A) \gpd M$. 
It is a Lie subgroupoid of the frame groupoid $\Phi(L)$ and it is called the {\em groupoid of inner automorphisms} of $L$. Let 
$F \rightarrow M$ be the Lie group bundle $L$ integrates to. 
Then the Lie algebroid coupling crossed module $\langle L,\ad,\ad(A),\rho_{*} \rangle$ integrates to a coupling crossed 
module of Lie groupoids, namely $\langle F,I,\Int(A),\rho \rangle$. Here $I : F \rightarrow \Int(A)$ maps an element $f$ which 
belongs to the fiber $F_{x}$ to the inner automorphism $I_{f}$ of the fiber $F_{x}$. The image of this map is the Lie group bundle 
$\Inn(F)$ of inner automorphisms of the fibers of $F$ and it is immediate that it differentiates to $\ad(L)$. Finally, the 
representation $\rho : \Int(A) * F \rightarrow F$ is $\rho(\phi,f) = \phi(f)$ for all $(\phi,f) \in \Int(A) * F$.
\begin{rmk}
If a Lie algebroid $A$ over $M$ integrates to a Lie groupoid $G$ then we may quantize the natural Poisson structure of $A^*$. An account of this is given in \cite[III, 3.11]{Landsman}. Very roughly, one considers the tangent groupoid $G_{T} = A \times \{ 0 \} \cup G \times \R^{*}$ over $M \times \R$ (where $A$ is considered to be a Lie groupoid with the fiberwise additive structure). This induces a natural short exact sequence of $C^*$-algebras $0 \to C^{\infty}(A^*) \to C^*(C_{T}) \to C^{*}(G) \otimes C_{0}(\R) \to 0$. It turns out that the commutators of $C^*(G)$ approximate asymptotically the Poisson bracket of $C^{\infty}(A^*)$.

Notice that the above method reduces the possibly non-integrable Lie algebroid $A$ to the integrable $\ad(A)$. Therefore the Poisson structure of $\ad(A)^*$ is always quantizable with the previous method, although the Poisson structure of $A^*$ may not be. Moreover, by dualising $\ad$ we get an injection $\ad^* : \ad(A)^{*} \to A^*$ which is a Poisson map. This can easily be extended to any regular Lie algebroid.
\end{rmk}

\begin{sloppypar}
\begin{prop}
The obstruction to the integrability of the Lie algebroid $A$ is $\Obs\langle F,I,\Int(A),\rho \rangle$.
\end{prop}
\end{sloppypar}

\pf If $\Obs\langle F,I,\Int(A),\rho \rangle = 0$, then the coupling crossed module of Lie groupoids has operator extensions which 
differentiate to operator extensions of the coupling crossed module of Lie algebroids $\langle L,\ad,\ad(A),\rho_{*} \rangle$. Therefore 
there exists a Lie groupoid which integrates the Lie algebroid $A$. \boom

The following result is an immediate consequence of the previous proposition and \ref{constant:opextobstr}.
\begin{cor}\label{constant:integrobstr}
The integrability obstruction of a transitive Lie algebroid takes values in \v{C}ech cohomology with constant coefficients. 
\end{cor}
 
\section{Crossed modules with PBG structures}

The rest of this paper intends to generalise the crossed module approach described above to general extensions of integrable transitive Lie algebroids by Lie algebra bundles. To this end, we start by recalling briefly the reformulation of such extensions given in \cite{Mackenzie:1988OEPB}, \cite{Mackenzie:1989IOELA} and \cite[I \S2.5, \S4.5]{KCHM:new}.

\begin{df}
Let $P(M,G)$ be a principal bundle. A transitive Lie algebroid $A$ over $P$ is called a {\em PBG-algebroid} if the Lie group $G$ acts on the manifold $A$ so that for every $g \in G$, the diagram
\begin{diagram}
A & \rTo^{R_{g}} & A \\
\dTo &  & \dTo \\
P & \rTo^{R_{g}} & P
\end{diagram}
is an automorphism of Lie algebroids.
\end{df}
A PBG-algebroid $A$ over the principal bundle $P(M,G)$ is denoted by $A \Rightarrow P(M,G)$. If $A$ and $A'$ are PBG-algebroids over the same principal bundle $P(M,G)$, a morphism of PBG-algebroids is a morphism of Lie algebroids $\psi : A \rightarrow A'$ satisfying $\psi(Xg) = \psi(X)g$ for all $X \in A$ and $g \in G$. A Lie algebra bundle $K \rightarrow P(M,G)$ is a PBG-Lie algebra bundle if, as a totally intransitive Lie algebroid, it is a PBG-algebroid. Given a PBG-Lie algebra bundle $K \rightarrow P(M,G)$, it is straightforward that $\CDO[K]$ is itself a PBG-algebroid. 
A {\em representation of PBG-algebroids} is a morphism of PBG-algebroids $\rho : A \rightarrow \CDO[K]$.

In \cite{Mackenzie:1989IOELA} it is shown that PBG-algebroids correspond to extensions of integrable Lie algebroids by Lie algebra bundles. To give an outline of this correspondence, let us start by considering a PBG-algebroid $A$ over the principal bundle $P(M,G)$. This can be written as an extension of Lie algebroids as $K \inj A \surj TP$, where $K$ is a PBG-Lie algebra bundle over $P(M,G)$. It is shown in \cite{Androulidakis:2004CELG} that the quotient space $\frac{A}{G}$ is always a manifold, therefore the above extension gives rise to an extension of Lie algebroids $\frac{K}{G} \inj \frac{A}{G} \surj \frac{TP}{G}$ over $M$.

On the other hand, consider a transitive Lie groupoid $\Omega \gpd M$, a Lie algebra bundle $K$ over $M$, and an extension of Lie algebroids $K \inj A \surj A\Omega$. Choose a basepoint in $M$ and let $P(M,G)$ denote the principal bundle corresponding to the groupoid $\Omega$. Then the pull-back of $A$ over the bundle projection $P \rightarrow M$ is a PBG-algebroid over $P(M,G)$.

Note that the right-splittings of the extension $K \inj A \surj A\Omega$ correspond to those splittings of the pullback PBG-algebroid which are equivariant with respect to the group action.
\begin{df}
Let $A \Rightarrow P(M,G)$ be a PBG-algebroid. An {\em isometablic} connection of $A$ is a vector bundle morphism $\gamma : TP \rightarrow A$ such that
\begin{enumerate}
\item $q \circ \gamma = id_{TP}$, where $q$ is the anchor of $A$;
\item $\gamma(Xg) = \gamma(X)g$ for all $X \in TP$ and $g \in G$.
\end{enumerate}
\end{df}
The respective notion for PBG-Lie algebra bundles is the following:
\begin{df}
Let $K \rightarrow P(M,G)$ be a PBG-Lie algebra bundle. An {\em isometablic} Koszul connection of $K$ is a vector bundle morphism $\nabla : TP \rightarrow \CDO(K)$ such that
\[
\nabla_{Xg}(Vg) = [\nabla_{X}(V)]g
\]
for all $X \in TP$, $V \in K$ and $g \in G$.
\end{df}

It is natural to postulate a notion of crossed module which is compatible with the PBG structure.
\begin{df}
A {\em crossed module of PBG-algebroids} over the principal bundle $P(M,G)$ is a crossed module of Lie algebroids $\langle K,\da,A,\rho \rangle$, where $K \rightarrow P(M,G)$ is a PBG-Lie algebra bundle, $A \Rightarrow P(M,G)$ is a PBG-algebroid, $\da : K \rightarrow A$ is a morphism of PBG-algebroids and $\rho : A \rightarrow CDO[K]$ is a representation of PBG-algebroids.
\end{df}

Let us now give a brief outline of the correspondence of pair crossed modules of PBG-algebroids with general crossed modules of integrable Lie algebroids. To start with this, consider a pair crossed module of PBG-algebroids $\langle K, \da, A, \rho \rangle$ over the principal bundle $P(M,G,p)$. Let $L \rightarrow P(M,G,p)$ be the kernel of the anchor of $A$, itself a PBG-Lie algebra bundle. It was shown in \cite{Androulidakis:2004CELG} that, because both $K$, $L$ and $A$ are PBG as Lie algebroids, the quotient manifolds $\frac{K}{G}$, $\frac{L}{G}$ and $\frac{A}{G}$ exist. Also, since the map $\da : K \rightarrow L$ is equivariant, it quotients to a morphism of vector bundles $\da^{/G} : \frac{K}{G} \rightarrow \frac{L}{G}$. This is a surjective morphism of Lie algebra bundles because $\da$ itself is surjective morphism of Lie algebra bundles. The representation $\rho : A \rightarrow \CDO[K]$ induces a representation $\rho^{/G} : \frac{A}{G} \rightarrow \CDO[\frac{K}{G}]$, defined by
\[
\rho^{/G}(\langle X \rangle)(\langle \mu \rangle) = \langle \rho(X)(\mu) \rangle
\]
for all $X \in A$ and $\mu \in K$. The fact that $\da^{/G}$ and $\rho^{/G}$ satisfy the properties of definition \ref{df:xm} follows from the fact that $\da$ and $\rho$ satisfy the respective properties for a pair crossed module of PBG-algebroids. Therefore we get the following crossed module of Lie algebroids:
\begin{diagram}
\ker{\da^{/G}} & & & & & & \\
\dTeXembed & & \rho^{/G} :\frac{A}{G} & \rTeXto & \CDO[\frac{K}{G}] & \\
 \frac{K}{G} & & & & & \\
\dTeXonto<{\da^{/G}} & & & & & \\
\im(\da^{/G}) & \rTeXembed & \frac{A}{G} & \rTeXonto^{\natural} & \frac{TP}{G} & \\
\end{diagram}

Conversely, suppose $\Omega \gpd M$ is a transitive Lie groupoid and $\langle K,\da,A,\rho \rangle$ is a crossed module of $A\Omega$. Choose a basepoint in $M$ and denote $P(M,G,p)$ the principal bundle corresponding to $\Omega$. We will show that this crossed module induces a pair crossed module of PBG-algebroids over $P(M,G)$. First, take the extension of Lie algebroids $Im(\tau) \inj A \stackrel{\pi}{\surj} A\Omega$. It was proven in \cite{Mackenzie:1989IOELA} that this extension induces a PBG-algebroid $A \actgpd p \Rightarrow P(M,G)$, where $A \actgpd p$ is the pullback vector bundle of $A$ over $p : P \rightarrow M$. That is $A \actgpd p = \{ (u,X) \in P \times A : X \in A_{p(u)} \}$. This Lie algebroid can be realised as an extension in the following way:
\[
(Im(\tau)) \actgpd p \inj A \actgpd p \stackrel{q^{!}}{\surj} TP.
\]
Recall that $\Gamma(A \actgpd p) \equiv C^{\infty}(P) \otimes \Gamma A$, where $f \otimes X$
corresponds to $f \cdot (X \circ p)$. The PBG-algebroid structure can be given both immediately
and on the module of sections. Thus the anchor is $q^{!}(u,X) = TR_{u}(\pi \circ X)$ and on the
sections, $q^{!}(f \otimes X)_{u} = f(u) \cdot TR_{u}(\pi(X_{p(u)}))$. The Lie bracket is
\[
[f \otimes X,h \otimes Y] = (f \cdot h) \otimes [X,Y] + (f \cdot \Vec{\pi(X)}(h)) \otimes Y -
(h \cdot \Vec{\pi(Y)}(f)) \otimes X,
\]
where $\Vec{\pi(X)}_{u} = TR_{u}(\pi(X)_{1_{\beta(u)}})$ is the right-invariant vector field on
$P$ corresponding to $\pi(X) \in \Gamma A\Omega$. We denote the action of $G$ on $A \actgpd p$ by
$R_{g}^{!}$ for every $g \in G$. This is defined as $R_{g}^{!}(u,X) = (u \cdot g,X)$ or, on the
sections by $R_{g}^{!}(f \otimes X) = (f \circ R_{g^{-1}}) \otimes X$.

The PBG-Lie algebra bundle $K \actgpd p$ is the pullback vector bundle of $K$ over $p: P
\rightarrow M$. Namely, $K \actgpd p = \{ (u,V) \in P \times K : V \in K_{p(u)} \}$. Again, from
the standard result for vector bundles we have $\Gamma(K \actgpd p) = C^{\infty}(P) \otimes_
{C^{\infty}(M)} \Gamma K$ where $f \otimes V$ corresponds to $f \cdot (V \circ p)$. 

The PBG-Lie algebra bundle structure of $K \actgpd p$ can also be given both immediately and on the sections. To show that $K \actgpd p$ is indeed a Lie algebra bundle, we have the following
theorem:

\begin{thm}
Suppose $K \rightarrow M$ is a vector bundle with a Lie bracket $[,]^{K} : M \rightarrow
Alt^{2}(K;K)$, $P$ a manifold and $p : P \rightarrow K$ a smooth map. Also, suppose $\nabla^{K} :
TM \rightarrow CDO[K]$ is a connection of $K$ such that
\[
\nabla^{K}_{X}[V,W] = [\nabla^{K}_{X}(V),W]^{K} + [V,\nabla^{K}_{X}(W)]^{K}.
\]
If $[,]^{!} P \rightarrow Alt^{2}(K \actgpd p, K \actgpd p)$
is the Lie bracket on $K \actgpd p$ defined as
\[
[(u,V),(u,W)]^{!} = (u,[V,W]^{K}_{p(u)}),
\]
then the map $\nabla^{!} : TP \rightarrow CDO[K \actgpd p]$ defined by $\nabla^{!}_{Y}(u,V) =
(u,\nabla^{K}_{Tp(Y)}(V))$ is a connection in $K \actgpd p$ and it satisfies
\[
\nabla^{!}_{Y}([(u,V),(u,W)]^{!}) = [\nabla^{!}_{Y}(u,V),(u,W)]^{!} +
[(u,V),\nabla^{!}_{Y}(u,W)]^{!}
\]
\end{thm}

\pf~It suffices to prove that $\nabla^{!}$ satisfies the last equality. This follows immediately
from the respective equality for $\nabla^{K}$. \boom
{\bf Remark.} The connection $\nabla^{!}$ of the previous theorem on the section-level is given
by the formula:
\[
\nabla^{!}_{Y}(f \otimes V) = f \otimes \nabla^{K}_{Tp(Y)}(V) + Tp(Y)(f) \otimes V.
\]

On the level of sections, the expression for the Lie bracket of $K \actgpd p$ defined above is given by $[f \otimes V,h \otimes W]^{!} = f \cdot h \otimes [V,W]$.

We now give an isometablic version of \cite[II 6.4.5]{KCHM:new}. It gives a criterion for the existence of a PBG-Lie algebra bundle structure on a vector bundle.
\begin{thm} \label{2.5.3}
Let $K$ be a vector bundle over $P(M,G)$ on which $G$ acts by isomorphisms and $[\ ,\ ]$ a section of the vector bundle $Alt^{2}(K;K)$. Then the following three conditions are equivalent:
\begin{enumerate}
\item The fibers of $K$ are pairwise isomorphic as Lie algebras.
\item $K$ admits an isometablic connection $\nabla$ such that 
\[
\nabla_{X}[V,W] = [\nabla_{X}(V),W] + [V,\nabla_{X}(W)] 
\]
for all $X \in \Gamma TP$ and $V,W \in \Gamma K$.
\item $K$ is a PBG-Lie algebra bundle.
\end{enumerate}
\end{thm}

\begin{cor}
If $K \rightarrow M$ is a Lie algebra bundle and $P(M,G,p)$ a principal bundle then $K \actgpd p$
is a PBG-Lie algebra bundle.
\end{cor}

\pf~Define the action of $G$ on $K \actgpd p$ to be $R_{g}^{!}(u,V) = (u \cdot g,V)$. On the
section-level, it will be $R^{!}_{g}(f \otimes V) = (f \circ R_{g^{-1}}) \otimes V$.
We showed in \ref{2.5.3} that a vector bundle is a PBG-Lie algebra bundle if and only if it has an
isometablic Lie connection. If $K \rightarrow M$ has a Lie connection $\nabla^{K}$
then the connection $\nabla^{!}$ constructed in the previous theorem is also a Lie connection.
It is moreover isometablic because:
\[
\nabla^{!}_{TR_{g}(Y)}(u \cdot g,V) = (u \cdot g, \nabla^{K}_{Tp \circ TR_{g}(Y)}(V)) =
(u \cdot g,\nabla^{K}_{T(p \circ R_{g})}(V)) = (u \cdot g,\nabla^{K}_{Tp(Y)(V)}).
\]
\boom

The next step is to define the morphism of PBG-algebroids $\tau^{!} : K \actgpd p \rightarrow
A \actgpd p$. This is defined by $\tau^{!}(u,V) = (u,\tau(V))$, or, on the section-level
by $\tau^{!}(f \otimes V) = f \otimes \tau(V)$. It is a straightforward calculation to show that it is a morphism of PBG-algebroids.

Finally we need to define a representation $\rho^{!} : A \actgpd p \rightarrow CDO[K \actgpd p]$
and show that it satisfies the properties of definition 4.1.1. This is defined as $\rho^{!}(u,X)
(u,V) = (u,\rho(X)(V))$, or, on the section-level as
\[
\rho^{!}(f \otimes X)(h \otimes V) = (f \cdot h) \otimes \rho(X)(V) + (f \cdot \Vec{\pi(X)}(h))
\otimes V.
\]
Again, the proof that $\da^{!}$ and $\rho^{!}$ satisfy the necessary properties which make $\langle K \actgpd p, \da^{!}, A\actgpd p, \rho^{!} \rangle$ a pair crossed module of PBG-algebroids is a straightforward calculation. These considerations can be formulated to the following result:
\begin{thm}
Pair crossed modules of PBG-algebroids are equivalent to crossed modules of integrable Lie algebroids.
\end{thm}

Therefore, it suffices to work with pair crossed modules of PBG-algebroids. As far as their operator extensions are concerned, these are pairs $(A,\mu)$ where $A \Rightarrow P(M,G)$ is a PBG-algebroid and $\mu$ is a morphism of PBG-algebroids.
\begin{rmk}
Following the same process as the one described in \cite[II\S7.3]{KCHM:new}, and working with isometablic connections instead, the obstruction to the existence of an operator extensions is an element of $G$-equivariant Lie algebroid cohomology ${\cal{H}}^{3}_{G}$ and if it vanishes the operator extensions are classified by ${\cal{H}}^{2}_{G}$.
\end{rmk}

\subsection{Crossed modules of PBG-groupoids}

Every extension $K \inj A \surj A\Omega$ of transitive Lie algebroids over $M$ (where $A\Omega$ integrates to the Lie groupoid $\Omega \gpd M$) gives rise to a Lie algebroid crossed module of $A\Omega$, and in the previous section we showed that such crossed modules correspond to pair crossed modules of PBG-algebroids. Therefore, bearing in mind the ideas explained in section 3, the obstruction to the integrability of a general extension $K \inj A \surj A\Omega$ must coincide with the obstruction associated with a certain crossed module on the groupoid level, which involves a PBG structure as well. 

In this section we give the definition of this particular crossed module and its operator extensions, and show that such crossed modules correspond to pair crossed modules of PBG-algebroids via the processes of differentiation and integration. Let us start with a brief account on the prerequisites of PBG structures on the groupoid level.

\begin{df}\label{df:PBGgpd}
A {\em PBG-groupoid} is a Lie groupoid $\Upsilon \gpd P$ whose base is the total
space of a principal bundle $P(M,G)$ together with a right action of $G$ on 
the manifold $\Upsilon$ such that for all $(\xi,\eta) \in \Upsilon \times \Upsilon$ such that $s\xi = t\eta$ and $g \in G$ we have:
\begin{enumerate}
\item $t(\xi \cdot g) = t(\xi) \cdot g$ and $s(\xi \cdot g) =
s(\xi) \cdot g$
\item $1_{u \cdot g} = 1_{u} \cdot g$
\item $(\xi \eta) \cdot g = (\xi \cdot g)(\eta \cdot g)$
\item $(\xi \cdot g)^{-1} = \xi^{-1} \cdot g$
\end{enumerate}
\end{df}
The properties of a PBG-groupoid imply that the right translation on $\Upsilon$ is a Lie groupoid automorphism over the right translation of the principal bundle. A morphism $\phi$ 
of Lie groupoids between two PBG-groupoids $\Upsilon$ and $\Upsilon'$ over the same 
principal bundle is called a morphism of PBG-groupoids if it preserves the 
group actions. Namely, if $\phi \circ \Tilde{R}_{g} = \Tilde{R}_{g}' \circ \phi$ 
for all $g \in G$. In the same fashion, a PBG-Lie group bundle (PBG-LGB) is a 
Lie group bundle $F$ over the total space $P$ of a principal bundle $P(M,G)$ 
such that the group $G$ acts on $F$ by Lie group bundle automorphisms. We 
denote a PBG-LGB by $F \rightarrow P(M,G)$. It is easy to see that the gauge 
group bundle $I\Upsilon$ of a PBG-groupoid $\Upsilon \gpd P(M,G)$ is a PBG-LGB. It is straightforward that PBG-groupoids differentiate to PBG-algebroids.

The class of transitive PBG-groupoids is of interest here, and that is because these groupoids are equivalent to extensions of transitive Lie groupoids. Namely, given a transitive PBG-groupoid $\Upsilon \gpd P(M,G)$, its corresponding extension $I\Upsilon \inj \Upsilon \surj P \times P$ can be quotiened by $G$ (see \cite{Mackenzie:1989IOELA}) to give rise to the extension of transitive Lie groupoids over $M$.
\[
\frac{I\Upsilon}{G} \inj \frac{\Upsilon}{G} \surj \frac{P \times P}{G}.
\]
On the other hand, given an extension of transitive Lie groupoids $F \inj \Omega \surj \Phi$ over $M$, choose a basepoint and consider the corresponding extension of principal bundles $N \inj Q(M,H) \surj P(M,G)$. This gives rise to the principal bundle $Q(P,N)$, and in turn this forms the Lie groupoid $\Upsilon = \frac{Q \times Q}{N} \gpd P$. Now the Lie group $G$ acts on $\Upsilon$ by
\[
\langle q_{2},q_{1} \rangle g = \langle q_{2}h,q_{1}h \rangle,
\]
where $h$ is any element of $H$ which projects to $g$. A detailed account of these constructions can be found in \cite{Mackenzie:1989IOELA}, as well as the proof that they are mutually inverse. 

\begin{df}\label{df:xmpbggpds}
\begin{sloppypar}
A {\em crossed module of PBG-groupoids} is a quadruple $(F,\da,\Omega,\rho)$, where $\Omega \gpd P(M,G)$ is a PBG-groupoid, $\pi : F \rightarrow P(M,G)$ is a PBG-Lie group bundle, $\da : F \rightarrow \Omega$ is a morphism of PBG-groupoids over $P(M,G)$ and $\rho$ is a representation of $\Omega$ on $F$, all such that
\end{sloppypar}
\begin{enumerate}
\item $\rho(\xi g, f g) = \rho(\xi,f)g$ for all $(\xi,f) \in \Omega * F$ and $g \in G$;
\item $\da(\rho(\xi,f)) = \xi \da(f) \xi^{-1}$ for all $(\xi, f) \in \Omega * F$;
\item $\rho(\da(f),f') = ff'f^{-1}$ for all $f, f' \in F$ with $\pi(f) = \pi(f')$;
\item $Im(\da)$ is a closed embedded submanifold of $\Omega$.
\end{enumerate}
\end{df}
In the same fashion, $Im(\da)$ is a PBG-Lie group bundle which lies entirely in $I\Omega$ and is normal in $\Omega$, and the cokernel $\frac{\Omega}{Im(\da)}$ is a PBG-groupoid over $P(M,G)$. If the cokernel of a crossed module of PBG-groupoids is the pair groupoid $P \times P$, then the crossed module is called {\em pair}. If, moreover, $\ker\da = ZF$, then it is called a coupling. In the remaining of this paper we will be concerned only with pair crossed modules of PBG-groupoids.
\begin{df}
An {\em operator extension} of a {\em pair} crossed module of PBG-groupoids $\langle F, \da, \Omega, \rho \rangle$ over the principal bundle $P(M,G)$ is a pair $(\Phi,\mu)$ such that $\Phi$ is a PBG-groupoid over $P(M,G)$, $\mu : \Phi \rightarrow \Omega$ is a morphism of PBG-groupoids, and the pair is an operator extension in the sense of \ref{df:opextlalgd}.
\end{df}

\subsection{Differentiation}

Now consider a pair crossed module of PBG-groupoids $pxm = (F,\tau,\Omega,\rho)$ over the
principal bundle $P(M,G)$. From definition \ref{df:xmpbggpds} we then have:
\begin{enumerate}
\item $\rho(\xi g,fg) = \rho(\xi,f)g$ for all $(\xi,f) \in \Omega * F$ and $g \in G;$
\item $\tau(\rho(\xi,f)) = \xi \cdot \tau(f) \cdot \xi^{-1}$ for all
$(\xi,f) \in \Omega * F;$
\item $\rho(\tau(f),f') = ff'f^{-1}$ for all $f,f' \in F$ with $\pi(f) = \pi(f')$.
\end{enumerate}
In order to differentiate $pxm$ to a pair crossed module of PBG-algebroids, consider the
PBG-Lie algebra bundle $F_{*} \rightarrow P(M,G)$, the PBG-algebroid $A\Omega \Rightarrow P(M,G)$
and the morphism of PBG-Lie algebra bundles $\tau_{*} : F_{*} \rightarrow L\Omega$.

First of all we construct a representation $\rho_{*} : A\Omega \rightarrow CDO[F_{*}]$ which
preserves the $G-$actions. Since $\rho$ is an equivariant representation we have that
$\rho(\xi) : F_{\alpha(\xi)} \rightarrow F_{\beta(\xi)}$ is a Lie group isomorphism for all
$\xi \in \Omega$ such that for every $g \in G$ the isomorphism $\rho(\xi g) : F_{\alpha(\xi g)}
\rightarrow F_{\beta(\xi g)}$ is equal to $\rho(\xi)g$. Applying the Lie functor we have that
$(\rho(\xi g))_{*} = (\rho(\xi))_{*} g$ for all $g \in G$. Thus, we get a well defined morphism
of PBG-groupoids
\[
\tilde{\rho} : \Omega \rightarrow \Pi[F_{*}],~\xi \mapsto (\rho(\xi))_{*}.
\]
Denote $\rho_{*} : A\Omega \rightarrow CDO[F_{*}]$ the morphism of PBG-algebroids $\Tilde{\rho}$
differentiates to. This is the representation we are looking for. It is straightforward to show that $\rho_{*}$ and $\da_{*}$ satisfy the properties of lemma \ref{lemma:diffxm} (same proof), thus making $(F_{*},\da_{*},A\Omega,\rho_{*})$ a pair crossed module of PBG-algebroids.

Next, suppose that $pxm$ has an operator extension $(F \stackrel{\iota}{\inj} \Phi
\stackrel{(\beta,\alpha)}{\surj} P \times P,\mu)$. That is to say that $\mu : \Phi \rightarrow
\Omega$ is a morphism of PBG-groupoids such that the diagram
$$
\begin{diagram}
F & \rTeXembed^{\iota} & \Phi & \rTeXonto^{(\beta,\alpha)} & P \times P \\
\dTeXonto<{\tau} & & \dTeXonto>{\mu} & & \dequal \\
I\Omega & \rTeXembed & \Omega & \rTeXonto & P \times P \\
\end{diagram}
$$
commutes and $(\iota \circ \rho \circ \mu)(\omega) = I_{\omega} \circ \iota$ for all
$\omega \in \Phi$. It is immediate that the diagram
$$
\begin{diagram}
F_{*} & \rTeXembed^{\iota_{*}} & A\Phi & \rTeXonto & TP \\
\dTeXonto<{\tau_{*}} & & \dTeXonto>{\mu_{*}} & & \dequal \\
L\Omega & \rTeXembed & A\Omega & \rTeXonto & TP \\
\end{diagram}
$$
commutes, and, in the same fashion as with $\tau_{*}(\rho_{*}(X)(V)) = [X,\tau_{*}(V)]$,
one can prove that
\[
(\iota_{*} \circ \rho_{*} \circ \mu_{*})(X')(V) = [X',\iota(V)]
\]
for all $X' \in A\Phi$ and $V \in F_{*}$.

\subsection{Integration}

Suppose given a pair crossed module of PBG-algebroids $pxm_{*} = (K,\tau_{*},A,\rho_{*})$
over the principal bundle $P(M,G)$. The general theory (\cite{KCHM:new}, \cite{Mackenzie:1989IOELA}) induces that
the PBG-Lie algebra bundle $K \Rightarrow P(M,G)$ integrates to a PBG-Lie group bundle $F \rightarrow P(M,G)$ with connected
and simply connected fibers. This section proves that if the PBG-algebroid $A$ integrates to a
PBG-groupoid $\Omega \gpd P(M,G)$ which is $\alpha$-connected and $\alpha$-simply connected, then the pair crossed module $pxm_{*}$ integrates to a pair crossed module of
PBG-groupoids $pxm = (F,\tau,\Omega,\rho)$ over the principal bundle $P(M,G)$.

Since the PBG-Lie algebra bundle $F \rightarrow P(M,G)$ has connected and simply connected fibers, \cite{Mackenzie-Xu} shows that $\tau_{*} : F_{*} \rightarrow L\Omega$ integrates uniquely to a morphism of PBG-groupoids $\tau : F \rightarrow I\Omega$ which is onto. Thus, all we need to
show in order to prove the integrability of $pxm_{*}$ is that $\rho_{*}$ integrates uniquely to
an equivariant representation $\rho : \Omega \rightarrow \Pi(F)$ such that:
\begin{enumerate}
\item $\rho(\tau(f)) = I_{f}$ for all $f \in F$ and
\item $(\tau \circ \rho)(\xi) = I_{\xi} \circ \tau$ for all $\xi \in \Omega$.
\end{enumerate}

Let us start with the integration of a representation of PBG-algebroids.
Consider a PBG-groupoid $\Omega \gpd P(M,G)$ which is $\alpha-$connected and
$\alpha-$simply connected, a PBG-Lie group bundle $F \rightarrow P(M,G)$ with simply connected fibers and an
equivariant representation $\rho_{*} : A\Omega \rightarrow CDO[F_{*}]$ of the PBG-algebroid
$A\Omega$ on the PBG-Lie algebra bundle $F_{*}$. Since $\Omega$ is supposed to be $\alpha-$connected and
$\alpha-$simply connected, \cite{Mackenzie-Xu} shows that this integrates uniquely to a morphism of
PBG-groupoids $\Bar{\rho} : \Omega \rightarrow \Pi(F_{*})$ such that $\Bar{\rho}_{*} = \rho_{*}$.
For every $\xi \in \Omega$ we then have an isomorphism of Lie algebras
\[
\Bar{\rho}(\xi) : (F_{*})_{\alpha(\xi)} \rightarrow (F_{*})_{\beta(\xi)}.
\]
Moreover, for all $g \in G$ and $\xi \in \Omega$ we have $\Bar{\rho}(\xi g) = \Bar{\rho}(\xi)g$
because $\rho_{*}$ is equivariant. Therefore, from the general Lie theory for every
$\xi \in \Omega$ there is a Lie group isomorphism $\rho(\xi) : F_{\alpha(\xi)} \rightarrow
F_{\beta(\xi)}$ such that $(\rho(\xi))_{*} = \Bar{\rho}(\xi)$. So we get a map $\rho : \Omega
\rightarrow \Pi(F)$ sich that $\rho(\xi g) = \rho(\xi)g$ for all $g \in G$.
\begin{prop}
The map $\rho$ is $C^{\infty}-$differentiable.
\end{prop}

\pf~We work locally to prove this. Let $\{ P_{i} \}_{i \in I}$ be an open cover of $P$. Then
$\Omega_{P_{i}}^{P_{i}} \cong P_{i} \times \tau(H) P_{i}$ and $\Pi(F)_{P_{i}} \cong
P_{i} \times Aut(H) \times P_{i}$, where $H$ is the fiber type of $F$. Then $\rho$ over $P_{i}$
is the map
\[
\rho(u,\tau(h),v) = (u,\theta(u) \circ f(\tau(h)) \circ \theta(v)^{-1},v)
\]
where $\theta : P_{i} \rightarrow Aut(H)$ is a map (not $C^{\infty}$) and $f : H \rightarrow
Aut(H)$ is a Lie group morphism. Also, $\Pi(F_{*})_{P_{i}} \cong P_{i} \times Aut(\hoh)
\times P_{i}$ anf $\Bar{\rho}$ over $P_{i}$ becomes
\[
\Bar{\rho}(u,h,v) = (u,\Bar{\theta}(u) \circ \Bar{f}(h) \circ \Bar{\theta}{v}^{-1},v)
\]
where $\Bar{\theta} : P_{i} \rightarrow Aut(\hoh)$ is a $C^{\infty}-$map and $\Bar{f} : H
\rightarrow Aut(\hoh)$ is a Lie group morphism. We know that $(\rho(u,h,v))_{*} = \Bar{\rho}
(u,h,v)$, therefore the map $P_{i} \times H \times P_{i} \rightarrow Aut(\hoh)$ defined by
$(u,h,v) \mapsto (\theta(u))_{*} \circ (f(h))_{*} \circ (\theta(v))_{*}^{-1}$ is smooth. It
follows that the map $P_{i} \times H \times P_{i} \rightarrow Aut(H)$ defined by
$(u,h,v) \mapsto \theta(u) \circ f(h) \circ \theta(v)^{-1}$ is smooth. That is because of a more
general result which says that a linear first order system of equations, whose right-hand sides
depend smoothly on auxiliary parameters, has solutions which depend smoothly on these parameters,
providing that the initial conditions vary smoothly. Therefore, $\rho$ is smooth.~\boom
\begin{prop}
The map $\rho$ is a morphism of Lie groupoids.
\end{prop}

\pf~All we need to prove is $\rho(\eta\xi) = \rho(\eta) \circ \rho(\xi)$ for all $(\eta,\xi) \in
\Omega * \Omega$. To this end, take $\eta,\xi \in \Omega$ such that $\alpha(\xi) = u$,
$\beta(\xi) = v = \alpha(\eta)$ and $\beta(\eta) = w$. Now consider the Lie group automorphism
\[
f_{\eta\xi} = \rho(\eta\xi) \circ (\rho(\xi))^{-1} \circ (\rho(\eta))^{-1} : F_{w} \rightarrow
F_{w}.
\]
We will show that $f_{\eta\xi} = id_{F_{w}}$. Indeed, if $e_{w}$ is the identity element in
$F_{w}$ we have:
\begin{multline*}
T_{e_{w}}F_{\eta\xi} = T_{e_{w}}(\rho(\eta\xi) \circ (\rho(\xi))^{-1} \circ (\rho(\eta))^{-1}) = \\
= T_{((\rho(\xi))^{-1} \circ (\rho(\eta))^{-1})(e_{w})} \rho(\eta\xi) \circ
T_{e_{w}}((\rho(\xi))^{-1} \circ (\rho(\eta))^{-1}) = \\
= T_{e_{w}}(\rho(\eta\xi)) \circ T_{(\rho(\eta))^{-1}(e_{w})}(\rho(\xi))^{-1} \circ
T_{e_{w}}(\rho(\eta))^{-1} = \\
= T_{e_{w}}(\rho(\eta\xi)) \circ [T_{e_{u}}(\rho(\xi))]^{-1} \circ [T_{e_{v}}(\rho(\eta))]^{-1} = \\
= \Bar{\rho}(\eta\xi) \circ (\Bar{\rho}(\xi))^{-1} \circ (\Bar{\rho}(\eta))^{-1} =
id_{(F_{*})_{w}}.
\end{multline*}
Since $F$ has connected and simply connected fibers we get $f_{\eta\xi} = id_{F_{w}}$, thus
$\rho$ is indeed a morphism of Lie groupoids.~\boom

Now we can proceed to the integration of the pair crossed module.
We need to prove that $\rho$ and $\theta$ satisfy the identities mentioned in the beginning. To
this end, we need to establish the PBG-Lie group bundle morphisms $I$ and $Ad$ and the PBG-Lie
algebra bundle morphism $ad$.

Consider the PBG-Lie group bundle $F \rightarrow P(M,G)$ with fiber type $H$ and let
$\{ \psi_{i} : P_{i} \times H \rightarrow F_{P_{i}} \}_{i \in I}$ be a section atlas of it. It is
easily verified that the Lie group bundle $Aut(F) \rightarrow P(M,G)$ is a PBG-Lie group bundle
and the family of maps
$\{ \psi_{i}^{Aut} : P_{i} \times Aut(H) \rightarrow Aut(F)_{P_{i}} \}_{i \in I}$ defined by
\[
\psi_{i}^{Aut}(u,\phi \in Aut(H)) = \psi_{i,u} \circ \phi \circ \psi_{i,u}^{-1}
\]
is a section atlas for this bundle.
\begin{prop}
The map $I : F \rightarrow Aut(F)$ defined by $I_{f}(f') = f f'f^{-1}$ for all $f,f' \in F$ such
that $\pi(f) = \pi(f')$ is a PBG-Lie group bundle morphism. Locally it is of the form $I^{i} : F_{P_{i}}
\rightarrow (Aut(F))_{P_{i}}$ where
\[
I^{i}(u,h) = (u,I^{H}_{h})
\]
for all $(u,h) \in P_{i} \times H \stackrel{\psi_{i}}{\cong} F_{P_{i}}$. Here $I^{H}$ is the
inner automorphism of $H$.
\end{prop}

\pf~Immediate.~\boom

Next we consider the PBG-Lie group bundle $Aut(F_{*}) \rightarrow P(M,G)$. The section atlas of this bundle is
$\{ (\psi_{i}^{Aut})_{*} : P_{i} \times Aut(\hoh) \rightarrow Aut(F_{*})_{P_{i}} \}_{i \in I}$
defined by
\[
(\psi_{i}^{Aut})_{*}(u,\phi_{*} \in Aut(\hoh)) = (\psi_{i,u})_{*} \circ \phi_{*} \circ
(\psi_{i,x}^{-1})_{*} = T_{e}(\psi_{i}^{Aut}(u,\phi))
\]
for all $i \in I,~u \in P_{i}$ and $\phi \in Aut(H)$.
\begin{prop}
The map $Ad : F \rightarrow Aut(F_{*})$ defined by $Ad_{f} = T_{e_{u}}I_{f}$ for all
$f \in F_{u},~u \in P$ is a PBG-Lie group bundle morphism. Locally it is of the form $Ad^{i} : F_{P_{i}}
\rightarrow Aut(F_{*})_{P_{i}}$ where
\[
Ad^{i}(u,h) = (u,Ad^{H}_{h})
\]
for all $(u,h) \in P_{i} \times H \stackrel{\psi_{i}}{\cong} F_{P_{i}}$. Here $Ad^{H}$ is the
adjoint representation on $H$.
\end{prop}

\pf~Immediate~\boom

The representation $Ad$ differentiates to the PBG-Lie algebra bundle morphism $ad : F_{*} \rightarrow
Der(F_{*})$ defined by
\[
ad_{V}(W) = [V,W] = (T_{e_{u}}Ad(V))(W)
\]
for all $V,W \in (F_{*})_{u},~u \in P$.

Now we can proceed to the proof of the first identity. For all $V,W \in F_{*}$ we have
$\rho_{*}(\tau_{*}(V))(W) = [V,W]$, or $(\Bar{\rho} \circ \tau)_{*} = Ad_{*}$. Since
$F$ has connected and simply connected fibers we have $\Bar{\rho}(\tau(f)) = Ad_{f}$ for all
$f \in F$. Therefore,
\[
\rho \circ \tau = I.
\]

For the second identity, take an $X \in \Gamma A\Omega$. Then $X$ induces a vector field
$\Vec{X} \in {\cal{X}}(\Omega)$ which is defined as $\Vec{X}_{\xi} = T_{1_{\beta(\xi)}}R_{\xi}
(X_{\beta(\xi)})$ for all $\xi \in \Omega$. This is an $\alpha-$vertical ($\Vec{X}_{\xi} \in
T_{\xi}\Omega_{\alpha(\xi)}$) and right-invariant ($\Vec{X} \circ R_{\xi} = TR_{\xi} \circ
\Vec{X}$) vector field on $\Omega$. Let $\phi : (-\epsilon,\epsilon) \times {\cal{U}}_{0}
\rightarrow {\cal{V}}_{0}$ be the flow of $\Vec{X}$, where ${\cal{U}}_{0},{\cal{V}}_{0} \subseteq
\Omega$. Then, it is immediate that every $\phi_{t} : {\cal{U}}_{0} \rightarrow {\cal{V}}_{0}$ has
the properties $\alpha \circ \phi_{t} = \alpha$ and $\phi_{t} \circ R_{\xi} = R_{\xi} \circ
\phi_{t}$ for all $\xi \in \Omega$ and $t \in (-\epsilon,\epsilon)$.

Denote ${\cal{U}} = \beta({\cal{U}}_{0})$ and ${\cal{V}} = \beta({\cal{V}}_{0})$ and let
$\Psi : (- \epsilon,\epsilon) \times {\cal{U}} \rightarrow {\cal{V}}$ be the map $\psi_{t}(u) =
\beta(\phi_{t}(\eta))$ for all $\eta \in {\cal{U}}_{0}^{u}$. This is well defined because if we
consider an $\eta' \in {\cal{U}}_{0}^{u}$ then there is a $\xi \in {\cal{U}}_{0}$ such that
$\eta' = \eta \cdot \xi$. Consequently,
\[
\psi_{t}(u) = \beta(\phi_{t}(\eta\xi)) = \beta(\phi_{t}(\eta) \cdot \xi) = \beta(\phi_{t}(\eta)).
\]
Finally, for all $t \in (-\epsilon,\epsilon)$ we have $\alpha \circ \phi_{t} = \alpha$, $\beta
\circ \phi_{t} = \psi_{t} \circ \beta$ and $\phi_{t}(\xi\eta) = \phi_{t}(\xi) \cdot \eta$.
Therefore \cite[I 1.4.12]{KCHM:new} shows that $\phi_{t}$ is the restriction
to ${\cal{U}}_{0}$ of a unique local left-translation $L_{\sigma_{t}} : \Omega^{\cal{U}}
\rightarrow \Omega^{\cal{V}}$, where
\[
\sigma_{t}(u) = \phi_{t}(\xi) \cdot \xi^{-1}
\]
for all $\xi \in {\cal{U}}_{0}^{u}$. We define the exponential map $Exp : (-\epsilon,\epsilon)
\times \Gamma A\Omega \rightarrow \Gamma_{\cal{U}}\Omega$ by
\[
ExptX = \sigma_{t}.
\]

Now take an $X \in \Gamma A\Omega$ and a $V \in \Gamma F_{*}$. From the properties of the
exponential, for all $u \in P$ we have:
\begin{multline*}
\tau_{*}(\rho(X)(V_{u})) = -\frac{d}{dt}\tau_{*}(\Bar{\rho}(ExptX(u))(V_{u})) \mid_{0} = \\
= -\frac{d}{dt}\tau_{*}(\Bar{\rho}(\phi_{t}(\xi) \cdot \xi^{-1})(V_{u}))\mid_{0} =
-\frac{d}{dt}\tau_{*}(\Bar{\rho}(\phi_{t}(\xi))[\Bar{\rho}(\xi^{-1})(V_{u})])\mid_{0} = \\
= -\frac{d}{dt}(\tau \circ \rho(\phi_{t}(\xi)))_{*}[\Bar{\rho}(\xi)^{-1}(V_{u})]\mid_{0}
\end{multline*}
and
\begin{multline*}
ad_{X}(\tau_{*}(V_{u})) = -\frac{d}{dt}Ad_{(ExptX(u))}(\tau_{*}(V_{u})) \mid_{0} = \\
= -\frac{d}{dt}(I_{\sigma_{t}(u)})_{*}(\tau_{*}(V_{u}))\mid_{0} =
-\frac{d}{dt}(I_{\phi_{t}(\xi)\xi^{-1}})_{*}(\tau_{*}(V_{u}))\mid_{0} = \\
= -\frac{d}{dt}(I_{\phi_{t}(\xi)})_{*}[(I_{\xi^{-1}} \circ \tau)_{*}(V_{u})]\mid_{0}.
\end{multline*}
Consider the curves $\delta_{u},\gamma_{u} : (-\epsilon,\epsilon) \rightarrow F_{*}$ defined by
\[
\gamma_{u}(t) = (\tau \circ \rho(\phi_{t}(\xi)))_{*}[\Bar{\rho}(\xi)^{-1}(V_{u})]
\]
and
\[
\delta_{u}(t) = (I_{\phi_{t}(\xi)})_{*}[(I_{\xi^{-1}} \circ \tau)_{*}(V_{u})].
\]
Obviously, $\gamma_{u}(0) = \delta_{u}(0) = \tau_{*}(V_{u})$. Since $\tau_{*}(\rho(X)(V))
= [X,\tau_{*}(V)]$ we have
\[
\frac{d}{dt}\gamma_{u}(t)\mid_{0} = \frac{d}{dt}\delta_{u}(t) \mid_{0}.
\]
\begin{lem}
There is a $\delta < \epsilon$ such that $\frac{d}{dt}\gamma_{u}(t)\mid_{t_{0}} =
\frac{d}{dt}\delta_{u}(t)\mid_{t_{0}}$
for all $|t_{0}|<\delta$.
\end{lem}

\pf~For all $t \in (-\epsilon,\epsilon)$ we have:
\begin{multline*}
\gamma_{u}(t) = (\tau \circ \rho(\phi_{(t-t_{0}) + t_{0}}(\xi)))_{*}[\Bar{\rho}(\xi)^{-1}
(V_{u})] 
= (\tau \circ \rho(\phi_{t-t_{0}}(\xi)))_{*}[\Bar{\rho}(\xi)^{-1}(V_{u})].
\end{multline*}
Therefore,
\begin{multline*}
\gamma_{u}(t) - \gamma_{u}(t_{0}) =
(\tau \circ \rho(\phi_{t-t_{0}}(\xi)))_{*}[\Bar{\rho}(\xi)^{-1}(V_{u})]  -
(\tau \circ \rho(\phi_{t_{0}}(\xi)))_{*}[\Bar{\rho}(\xi)^{-1}(V_{u})] = \\
= (\tau \circ \rho(\phi_{t-t_{0}}(\xi)))_{*}[\Bar{\rho}(\xi)^{-1}(V_{u})] =
\gamma_{u}(t-t_{0}).
\end{multline*}
And of course, the same is true for $\delta_{u}$. Define $\Tilde{\gamma}_{u}(t) =
\gamma_{u}(t - t_{0})$ and $\Tilde{\delta}_{u}(t) = \delta_{u}(t-t_{0})$. Then,
\begin{multline*}
\frac{d}{dt}\gamma_{u}(t)\mid_{t_{0}} = lim_{t \rightarrow t_{0}}\frac{\gamma_{u}(t) -
\gamma_{u}(t_{0})}{t - t_{0}} = lim_{t\rightarrow 0}\frac{\Tilde{\gamma}_{u}(t)}{t} = \\
= \frac{d}{dt}\Tilde{\gamma}_{u}(t)\mid_{0} = \frac{d}{dt}\tilde{\delta_{u}(t)}\mid_{0} = ... =
\frac{d}{dt}\delta_{u}(t)\mid_{t_{0}}.
\end{multline*}
\boom

So, for all $|t|<\delta<\epsilon$ we have $(\tau \circ \rho(ExptX))_{*} =
(I_{ExptX} \circ \tau)_{*}$ and since $\Omega$ is $\alpha-$connected and $\alpha-$simply
connected we finally get
\[
\tau \circ \rho(ExptX) = I_{ExptX} \circ \tau
\]
for all $|t| < \delta$. Hence the desired equality. Combining this result with Proposition 4.4
of \cite{Mackenzie:1989IOELA} we get the following result.
\begin{thm}
Suppose given a PBG-Lie group bundle $F$ and a PBG-groupoid $\Omega$, both over the same 
principal bundle $P(M,G)$. Then any pair crossed module of PBG-algebroids 
$(F_{*},\tau_{*},A\Omega,\rho_{*})$ integrates to a pair crossed module of PBG-groupoids 
$(F,\tau,\Omega,\rho)$.
\end{thm}

\section{Classification of PBG-groupoids}

In order to characterize cohomologically the obstruction associated with a pair crossed module of PBG-groupoids, it is necessary to classify such groupoids. The classification given in \cite{Androulidakis:2004CELG} will be used in Section 7, for the enumeration of operator extensions when the lifting obstruction vanishes. As we discussed in the introduction, in this section we give a different cohomological classification of PBG-groupoids, which is consistent with the classification of transitive Lie algebroids given in \cite[II\S8.2]{KCHM:new}.

Let $\Omega \gpd P(M,G)$ be a PBG-groupoid and $\{ P_{i} \equiv U_{i} \times G \}_{i \in I}$ an atlas of its base principal bundle. It was shown in \cite{Androulidakis:2004CELG} that for every $i \in I$ there exists a flat isometablic connection $\gamma_{i} : TP_{i} \rightarrow A\Omega_{P_{i}}$. More than that, it was shown that as a morphism of PBG-algebroids, every $\gamma_{i}$ integrates to a morphism of PBG-groupoids $\theta_{i} : P_{i} \times P_{i} \rightarrow \Omega_{P_{i}}^{P_{i}}$. Now fix a $u_{0} \in P$ and denote $H = \Omega_{u_{0}}^{u_{0}}$. For every $i \in I$ choose a $u_{i} \in P_{i}$ and an arrow $\xi_{i} \in \Omega_{u_{0}}^{u_{i}}$. Now define the maps 
\[
\sigma_{i} : P_{i} \rightarrow \Omega_{u_{0}},~ \sigma_{i}(u) = \theta_{i}(u,u_{i}) \cdot \xi_{i}.
\]
These are sections of $\Omega$ and they respect the $G$-action in the following sense:
\[
\sigma_{i}(ug) = [\sigma_{i}(u)g] \cdot (\xi_{i}^{-1}g) \cdot \sigma_{i}(u_{i}g).
\]
These sections give rise to a family of representations $\{ \phi_{i} : G \rightarrow Aut(H) \}_{i \in I}$ of $G$ on $H$, namely
\[
\phi_{i}(g)(h) = \sigma_{i}(u_{i}g)^{-1} \cdot (\xi_{i}g) \cdot (hg) \cdot (\xi_{i}g)^{-1} \cdot \sigma_{i}(u_{i}g).
\]
It was shown in \cite{Androulidakis:2004CELG} that these representations are local expressions of the automorphism action of $G$ on the Lie group bundle $I\Omega$.

If we begin with a different local family $\{\gamma_{i}'\}_{i \in I}$ of flat isometablic connections, there exist 1-forms $\ell_{i}^{*} : TP_{i} \rightarrow P_{i} \times \hoh_{i}$ such that $\gamma_{i}' = \gamma_{i} + \ell_{i}^{*}$. Here $\hoh_{i}$ is the Lie algebra of the Lie group $\Omega_{u_{i}}^{u_{i}}$. Therefore the $\ell_{i}^{*}$s integrate to maps $\ell_{i} : P_{i} \times P_{i} \rightarrow \Omega_{u_{i}}^{u_{i}}$ such that $\theta_{i}' = \theta_{i} + \ell_{i}$. Define
\[
r_{i} : P_{i} \rightarrow H,~ r_{i}(u) = \xi_{i}^{-1}\cdot \ell_{i}(u,u_{i}) \cdot \xi_{i}.
\]
Now the respective sections are related by $\sigma_{i}' = \sigma_{i} \cdot r_{i}$, and with respect to the $G$-action the $r_{i}$s satisfy
\[
r_{i}(ug) = \phi_{i}(g)(r_{i}(u)) \cdot r_{i}(u_{i}g).
\]
Last, the representations arising from $\sigma_{i}'$ and $\sigma_{i}$ are related by
\[
\phi_{i}'(g)(h) = r_{i}(u_{i}g)^{-1} \cdot \phi_{i}(g)(h) \cdot r_{i}(u_{i}g).
\]

Now, instead of classifying $\Omega$ by the transition functions $s_{ij} : P_{ij} \rightarrow H$ associated with the sections $\sigma_{i}$, let us consider the following maps:
\[
\chi_{ij} : P_{ij} \times P_{ij} \rightarrow H,~ \chi_{ij}(u,v) = s_{ij}(u) \cdot s_{ji}(v)
\]
and
\[
\alpha_{ij} : P_{ij} \rightarrow Aut(H),~ \alpha_{ij}(u)(h) = s_{ij}(u) \cdot h \cdot s_{ji}(u).
\]
The $\alpha_{ij}$s are the transition functions of the PBG-Lie group bundle $I\Omega$. Together with the $\chi_{ij}$s they satisfy:
\begin{enumerate}
\item $\chi_{ik}(u,v) = \chi_{ij}(u,v) \cdot \alpha_{ij}(v)(\chi_{jk}(u,v))$.
\item For a choice of $u_{ij} \in P_{ij}$, $\alpha_{ij}(u) = I_{\chi_{ij}(u,u_{ij})} \circ I_{s_{ij}(u_{ij})}$.
\item $\chi_{ij}(ug,vg) = \phi_{i}(g)(\chi_{ij}(u,v))$.
\item $\alpha_{ij}(ug)(\phi_{j}(g)(h)) = \phi_{i}(g)(\alpha_{ij}(u)(h))$.
\end{enumerate}

\begin{df}\label{df:isomstd}
A pair $(\chi,\alpha)$ satisfying (i)--(iv) is called a {\em $\phi$-isometablic pair of transition data}.
\end{df}
The relation between two isometablic systems of transition data given in the following proposition was proven in \cite{Androulidakis:2004CELG}.
\begin{prop}\label{std:rel}
Two ${\cal{\phi}}$-isometablic and ${\cal{\phi'}}$-isometablic systems of transition data $({\cal{\chi}},{\cal{\alpha}})$ and $({\cal{\chi'}},{\cal{\alpha'}})$ respectively are related by
\begin{eqnarray}\label{systrelchi}
{\chi'}_{ij}(u,v) = r_{i}(u)^{-1}[\chi_{ij}(u,v) \cdot \alpha_{ij}(v)(r_{i}(u) \cdot r_{j}(v)^{-1})] \cdot r_{i}(v)
\end{eqnarray}
and
\begin{eqnarray}\label{systrelalpha}
{\alpha'}_{ij}(u) = I_{r_{i}(u)^{-1}} \circ \alpha_{ij}(u) \circ I_{r_{j}(u)}
\end{eqnarray}
\end{prop}
It is straightforward that the relation between isometablic systems of transition data is an equivalence relation, therefore it is legitimate to give the following definition.
\begin{df}
Two isometablic systems of transition data which satisfy ({\ref{systrelchi}}) and (\ref{systrelalpha}) are called {\em equivalent}.
\end{df}
Now we can proceed to show that isometablic systems of transition data classify PBG-groupoids.
\begin{prop}\label{PBGclass}
Suppose $P(M,G)$ be a principal bundle and $\{ U_{i} \}_{i \in I}$ is a simple open cover of $M$, whereas $\{ P_{i} \equiv U_{i} \times G \}_{i \in I}$ is an atlas over this cover. Let $\phi = \{ \phi_{i} : G \rightarrow Aut(H) \}_{i \in I}$ be a family of representations of $G$ on a Lie group $H$ and $(\chi,\alpha)$ a family of $\phi$-isometablic transition data. Then there exists a PBG-groupoid over $P(M,G)$ and a local family of flat isometablic connections which give rise to this data.
\end{prop}

\pf~For each $i \in I$, let $\Upsilon^{i} = P_{i} \times H \times P_{i}$ and on the disjoint sum of the $\Upsilon^{i}$s define an equivalence relation $\sim$ by
\[
(i,u,h,v) \cong (j,u',h',v') \Leftrightarrow u = u',~v = v',~h' = \chi_{ji}(u,v)\cdot\alpha_{ji}(v)(h).
\]
Denote the quotient set by $\Upsilon$ and the equivalence classes by $\langle i,(u,h,v) \rangle$. Define maps $\alpha, \beta : \Upsilon \rightarrow P$ by $\langle i, (u,h,v) \rangle \mapsto v$ and $\langle i, (u,h,v) \rangle \mapsto u$ respectively. The object inclusion map is $P_{i} \ni u \mapsto \langle i, (u,e_{H},u) \rangle$. it is easy to see that the map
\[
\Bar{\Psi}_{i} : P_{i} \times H \times P_{i} \rightarrow (\beta,\alpha)^{-1}(P_{i}),~(u,h,v) \mapsto \langle i, (u,h,v) \rangle
\]
is a bijection. Give $\Upsilon$ the smooth structure induced from the manifolds $P_{i} \times H \times P_{i}$ via $\Bar{\Psi}_{i}$.

Now we define a multiplication in $\Upsilon$. For $\xi, \eta \in \Upsilon$ such that $\alpha(\xi = \beta(\eta)) = u$, choose a $P_{i}$ containing $u$ and write $\xi = \langle i, (v,h,u) \rangle$, $\eta = \langle i, (u,h',w) \rangle$. Define
\[
\xi \cdot \eta = \langle i, (v,hh',w) \rangle.
\]
Finally, $G$ acts on $\Upsilon$ by
\[
\langle i, (u,h,v) \rangle g = \langle i, (ug,\phi_{i}(g)(h),vg) \rangle.
\]
It is left to the reader to verify that $\Upsilon$ is a well defined PBG-groupoid over $P(M,G)$. The PBG-algebroid it differentiates to is the one given in \cite[II\S5.4]{KCHM:new}, and the connections associated with the transition data we began with are the ones given there. \boom

{\bf Remark.} Note that for the previous construction the only property that we use is the cocycle condition that $(\chi, \alpha)$ satisfy, namely $\chi_{ij}(u,v) = \chi_{ik}(u,v) \cdot \alpha_{ik}(v)(\chi_{kj}(u,v))$, and that is to show that the relation $\sim$ is indeed an equivalence relation. On the other hand, the compatibility condition is not used here.
\\ \\ 
The following proposition shows that the PBG-groupoid arising from isometablic transition data is well defined up to equivalence. Its proof is a straightforward calculation.
\begin{prop}
Let $P(M,G)$ be a principal bundle, $\{ P_{i} \}_{i \in I}$ an open cover of $P$ by principal bundle charts, $H$ a Lie group and $\phi', \phi$ be two families of representations of $G$ on $H$ by 
which are equivalent under a family of maps $r = \{ r_{i} : P_{i} \rightarrow H \}_{i \in I}$
such that $r_{i}(ug) = \phi_{i}(g)(r_{i}(u)) \cdot r_{i}(u_{i}g)$ for all $u \in P_{i}, g \in G$ and $i \in I$. Let $(\chi',\alpha')$ and $(\chi,\alpha)$ be $\rho'$-isometablic and $\rho$-isometablic systems of transition data with values in $H$ respectively which are equivalent under the family of maps $r$. Let $\Omega'$ and $\Omega$ be the associated PBG-groupoids respectively. Then the map $\Xi : \Omega' \rightarrow \Omega$
defined by
\[
\langle i, (u, h, v) \rangle \mapsto \langle i, (u, r_{i}(u)^{-1} \cdot h \cdot r_{j}(v), v) \rangle
\]
is an isomorphism of PBG-groupoids over $P(M,G)$.
\end{prop}

\section{The obstruction of a pair crossed module of PBG-groupoids}

In this section we give the cohomological obstruction to the existence of an operator extension for a pair crossed module of PBG-groupoids. Let us start with such a crossed module $\langle F, \da, \Omega, \rho \rangle$ over the principal bundle $P(M,G)$. Then the PBG-groupoid $\Omega \gpd P(M,G)$ is the extension of PBG-groupoids
\[
Im(\da) \inj \Omega \stackrel{(\beta,\alpha)}{\surj} P \times P.
\]
Choose a simple open cover $\{ U_{i} \}_{i \in I}$ of $M$ and an atlas $\{ P_{i} \equiv U_{i} \times G \}_{i \in I}$ of the principal bundle, and consider a $\phi$-isometablic system of transition data $(\chi,\alpha)$. Note that in this context we denote $H$ the fiber type of $F$, therefore every $\phi_{i}$ is a representation of $G$ on $\da(H)$, namely $\phi_{i} : G \rightarrow Aut(\da(H))$. Now the following proposition shows that there exist canonical lifts of the representations $\phi_{i}$. Its proof is a straightforward calculation.
\begin{prop}
For every $i \in I$ the map $\Hat{\phi}_{i} : G \rightarrow Aut(H)$ defined by
\[
\Hat{\phi}_{i}(g)(h) = \rho(\sigma_{i}(u_{i}g)^{-1}\cdot (\xi_{i}g), hg)
\]
for all $g \in G$ and $h \in H$ is a representation of $G$ on $H$ and $\da \circ \Hat{\phi}_{i} = \phi_{i}$.
\end{prop}
The next two results show that there also exist canonical lifts of the transition functions $\alpha_{ij}$ of $Im(\da)$, to transition functions of $F$.

\begin{prop} \label{lift:PBGcharts}
Let $(F,\tau,\Omega,\rho)$ be a pair crossed module of PBG-groupoids. With the previous
notation, the maps $\psi_{i} : P_{i} \times H \rightarrow F_{P_{i}}$ defined by
$\psi_{i}(u,h) = \rho(\sigma_{i}(u),h)$ are charts of the Lie group bundle $F$ and they
are isometablic in the sense
\[
\psi_{i}(ug,\Hat{\phi}_{i}(g^{-1})(h)) = \psi_{i}(u,h) \cdot g
\]
for all $g \in G,~u \in P_{i}$ and $h \in H$.
\end{prop}

\pf~First of all, to prove that the $\psi_{i}$s are well defined, we need to ensure that the
restriction of $\rho$ on $\Omega_{u_{0}}^{P_{i}} * H$ takes values in $\pi^{-1}(P_{i})$.
Indeed, if $\xi \in \Omega_{u_{0}}^{P_{i}}$ and $f \in H$ is such that $\pi(f) = u_{0}$ then
$\pi(\rho(\xi,f)) = \beta(\xi) \in P_{i}$. The $\psi_{i}$s are injective because for all
$u,u' \in P_{i}$ and $f,f' \in H$ we have:
\begin{multline*}
\rho(\sigma_{i}(u),f) = \rho(\sigma_{i}(u'),f') \Rightarrow
\pi(\rho(\sigma_{i}(u),f)) = \pi(\rho(\sigma_{i}(u'),f')) \Rightarrow \\
\Rightarrow \beta(\sigma_{i}(u)) = \beta(\sigma_{i}(u')) \Rightarrow u = u'.
\end{multline*}
Since $\rho(\sigma_{i}(u))$ is an isomorphism on the fibers of $F$, we also have $f=f'$.

For the surjectivity of the $\psi_{i}$s, consider an $f \in F_{u} \subseteq F_{i}$ for some
$u \in P_{i}$. Then, because $\rho(\sigma_{i}(\pi(f)))$ is an isomorphism
$H \rightarrow F_{u}$, there is an $f' \in H$ such that $f = \rho(\sigma_{i}(\pi(f)),f')
= \psi_{i}(\pi(f),f')$. Last, the following diagram commutes
$$
\begin{diagram}
P_{i} \times N & & \rto^{\psi_{i}} & & \pi^{-1}(P_{i}) \\
               & \rdto<{pr_{1}} & & \ldto>{\pi} &      \\
               &      &    P_{i}   &       &      \\
\end{diagram}
$$
because $\pi(\psi_{i}(u,f)) = \pi(\rho(\sigma_{i}(u),f)) = \beta(\sigma_{i}(u)) = u$. For the isometablicity of the $\psi_{i}$'s we have:
\begin{multline*}
\psi_{i}(ug,\Hat{\phi}_{i}(g^{-1})(h)) = \rho(\sigma_{i}(ug),\Hat{\phi}_{i}(g^{-1})(h))
= \\ = \rho(\sigma_{i}(ug),\rho(\sigma_{i}(u_{i}g)^{-1} \cdot (\xi_{i}g),hg)) 
= \rho([\sigma_{i}(u)g] \cdot (\xi_{i}^{-1}g) \cdot \sigma_{i}(u_{i}g) \cdot
\sigma_{i}(u_{i}g)^{-1} \cdot (\xi_{i}g), hg) = \\
= \rho(\sigma_{i}(u)g,hg) = \rho(\sigma_{i}(u),h) \cdot g = \psi_{i}(u,h) \cdot g.
\end{multline*}
\boom
Now let us look at the transition functions of the Lie group bundle charts defined in the
previous theorem.
\begin{prop}
The transition functions of the charts $\{ \psi_{i} \}_{i \in I}$ are lifts of the transition
functions $\{ \Hat{\alpha}_{ij} \}_{i,j \in I}$, form a \v{C}ech-1-cocycle and are isometablic
with respect to the representations $\{ \Hat{\phi}_{i} \}_{i \in I}$.
\end{prop}

\pf~For all $u \in P_{ij}$ and $h \in H$, we have:
\begin{multline*}
\psi_{ij}(u)(h) = \psi_{i,u}^{-1}(\psi_{j,u}(h)) = \psi_{i,u}^{-1}
(\rho(\sigma_{j}(u),h))
= \rho(\sigma_{i}(u)^{-1} \cdot \sigma_{j}(u),h) 
= \rho(s_{ij}(u),h).
\end{multline*}
Therefore, $\tau(\psi_{ij}(u)(h)) = I_{s_{ij}(u)}(h) =
\alpha_{ij}(u)(h)$, so the $\psi_{ij}$'s are lifts of the $\alpha_{ij}$'s.
They form a \v{C}ech-1-cocycle because:
\begin{multline*}
[\psi_{jk}(u) \circ \psi_{ik}(u)^{-1} \circ \psi_{ij}(u)](h) =
\rho(s_{jk}(u),\rho(s_{ik}(u)^{-1},\rho(s_{ij}(u),h))) = \\
= \rho(s_{jk}(u) \cdot s_{ik}(u)^{-1} \cdot s_{ij}(u),h) = \rho(1_{u},h) = h
\end{multline*}
Moreover, they are isometablic with respect to the lifts $\{ \Hat{\phi}_{i} \}_{i \in I}$ of
the representations $\{ \phi \}_{i} \}_{i \in I}$ because:
\begin{multline*}
\psi_{ij}(ug)(\Hat{\phi}_{j}(g^{-1})(h)) = \rho(s_{ij}(ug),
\rho(\sigma_{j}(u_{j}g)^{-1} \cdot (\xi_{j}g),hg)) = \\
= \rho(\sigma_{i}(u_{i}g)^{-1} \cdot (\xi_{i}g) \cdot (s_{ij}(u)g) \cdot
(\xi_{j}^{-1}g) \cdot \sigma_{j}(u_{j}g) \cdot \sigma_{j}(u_{j}g)^{-1} \cdot
(\xi_{j}g), hg) = \\
= \Hat{\phi}_{i}(g^{-1})(\psi_{ij}(u)(h)).
\end{multline*}
\boom
Now let us show that there also exist canonical $\Hat{\phi}$-isometablic lifts of the $\chi_{ij}$s. Consider the quotient maps $\chi_{ij}^{/G} : \frac{P_{ij} \times P_{ij}}{G} \rightarrow \frac{\da(H)}{G}$ and $\da^{/G} : \frac{F}{G} \rightarrow \frac{Im(\da)}{G}$. The restriction of $\da^{/G}$ to $H = F_{u_{0}}$ is a a map $\da^{/G}\mid_{H} : \frac{H}{G} \rightarrow \frac{\da(H)}{G}$, where the action of $G$ on $\da(H)$ implied is $\phi_{i}$, and the action of $G$ on $H$ is $\Hat{\phi_{i}}$. This happens because the $\phi_{i}$s are local expressions of the $G$-action on $Im(F)$ as was shown in \cite{Androulidakis:2004CELG}. Bearing in mind that $P_{ij} \equiv U_{ij} \times G$, the quotient $\frac{P_{ij} \times P_{ij}}{G}$ is just $U_{ij} \times U_{ij}$, therefore we have the following diagram:
\begin{diagram}
 & & \frac{H}{G} \\
 & & \dto>{\da^{/G}} \\
U_{ij} \times U_{ij} & \rto_{\chi_{ij}^{/G}} & \frac{\da(H)}{G}
\end{diagram}
Note that $\frac{H}{G}(\frac{\da(H)}{G},\ker(\da), \da^{/G})$ is a principal bundle in a trivial way. Since the $U_{ij}$s are simply connected, it follows from \cite{Hu} that there exists a differentiable map $\Hat{\chi}_{ij}^{/G} : U_{ij} \times U_{ij} \rightarrow \frac{H}{G}$ such that the above diagram commutes. 

Denote $\sharp : P_{ij} \times P_{ij} \rightarrow U_{ij} \times U_{ij}$ and $\sharp^{H} : H \rightarrow \frac{H}{G}$ the natural projections. Since $\sharp^{H}$ is a pullback over the projection $\da^{/G}$ of the principal bundle $\frac{H}{G}(\frac{\da(H)}{G},\ker(\da), \da^{/G})$, there is a unique map $\Hat{\chi}_{ij} : P_{ij} \times P_{ij} \rightarrow H$ such that
\[
\sharp^{H} \circ \Hat{\chi}_{ij} = \Hat{\chi}_{ij}^{/G} \circ \sharp.
\]
Due to the $G$-invariance of $\sharp$ and $\sharp^{H}$, the map $\Hat{\phi}_{i}(g)^{-1} \circ \Hat{\chi}_{ij} \circ (R_{g} \times R_{g})$ also satisfies the previous equation for every $g \in G$, therefore it follows from the uniqueness argument that $\Hat{\chi}_{ij}$ is $\phi_{i}$-isometablic. these considerations consist the proof of the following result.
\begin{thm}
Let $\langle F, \da, \Omega, \rho \rangle$ be a pair crossed module of PBG-groupoids over a principal bundle $P(M,G)$ and $\phi = \{ \phi_{i} : G \rightarrow Aut(\da(H)) \}_{i \in I}$ a family of representations of $G$ on the image by $\da$ of the fiber type $H$ of the PBG-Lie group bundle $F$. Then there exists a canonical family of representations $\Hat{\phi} = \{ \Hat{\phi}_{i} : G \rightarrow Aut(H) \}_{i \in I}$ such that 
\begin{enumerate}
\item $\da \circ \Hat{\phi}_{i} = \phi_{i}$ for all $i \in I$;
\item For every $\phi$-isometablic system of transition data $(\chi,\alpha)$ of $\Omega$, there exists a canonical pair $(\Hat{\chi},\Hat{\alpha})$ with values in $H$, such that $\Hat{\alpha}$ is an isometablic cocycle of transition functions for the PBG-Lie group bundle $F$, and $\Hat{\chi}$ is a $\Hat{\phi}$-isometablic family of differential maps $\{ \Hat{\chi}_{ij} : P_{ij} \times P_{ij} \rightarrow H \}_{i,j \in I}$ such that $\da(\Hat{\chi},\Hat{\alpha}) = (\chi,\alpha)$.
\end{enumerate}
\end{thm}

If this lift of the transition data of $\Omega$ is a $\Hat{\phi}$-isometablic system of transition data itself, then it gives rise to a PBG-groupoid, with adjoint bundle $F$, and this would play the role of an operator extension for the given pair crossed module. We saw that this lift is indeed $\Hat{\phi}$-isometablic. As we remarked in \ref{PBGclass}, the only the only thing that is required is for the pair $(\Hat{\chi},\Hat{\alpha})$ to satisfy the cocycle condition. This can be reformulated to
\[
\psi_{i,v}(\Hat{\chi}_{ij}(u,v)) = \psi_{i,v}(\Hat{\chi}_{ik}(u,v)) \cdot \psi_{k.v}(\Hat{\chi}_{kj}(u,v)).
\]
Thus the failure of the pair $(\Hat{\chi},\Hat{\alpha})$ to satisfy the cocycle condition is the map $e_{ijk} : P_{ijk} \times P_{ijk} \rightarrow \ker\da \leq ZF$, defined by
\[
e_{ijk}(u,v) = \psi_{i,v}(\Hat{\chi}_{ij}(u,v)) \cdot [\psi_{k.v}(\Hat{\chi}_{kj}(u,v))]^{-1} \cdot \psi_{i,v}(\Hat{\chi}_{ik}(u,v)).
\]
The fact that it takes values in $\ker\da$ follows from the fact that the original system of transition data $(\chi,\alpha)$ does satisfy the cocycle condition. A routine calculation shows that $e_{ijk}(ug,vg) = e_{ijk}(u,v)g$ for all $g \in G$, and for $P_{ijkl} \neq \emptyset$ 
\[
e_{jkl} - e_{ikl} - e_{ijl} - e_{ijk} = 0 \in ZF
\]
and so $e$ is a 2-cocycle in $\check{H}^{2}_{G}(P \times P,ZF)$, the $G$-isometablic \v{C}ech cohomology of $P \times P$ with respect to the atlas $\{ P_{i} \cong U_{i} \times G \}_{i \in I}$ of the principal bundle $P(M,G)$, and with coefficients in the sheaf of germs of local isometablic maps from $P \times P$ to $ZF$. 

It is trivial to see that if a second family of lifts $P_{ij} \times P_{ij} \rightarrow H$ of the $\chi_{ij}$s is chosen then the resulting cocycle is cohomologous to $e$. More generally, if $(\chi',\alpha')$ is a second $\phi'$-isometablic system of transition data for $\Omega$, over the same atlas $\{ P_{i} \equiv U_{i} \times G \}_{i \in I}$ of the principal bundle $P(M,G)$, then it follows from the relations we gave in \ref{std:rel} and \ref{lift:PBGcharts} that $e'_{ijk} = e_{ijk}$.
\begin{thm}\label{thm:obstr}
Continuing the above notation, there exists an operator extension $(\Upsilon,\mu)$ for the pair crossed module of PBG-groupoids $\langle F,\da,\Omega,\rho \rangle$ iff $e = 0 \in \check{H}^{2}_{G}(P \times P,ZF)$. 
\end{thm}

\pf Assume that $e = 0$ and consider the PBG-groupoid $\Upsilon \gpd P(M,G)$ constructed directly from the pair $(\Hat{\chi},\Hat{\alpha})$ as in \ref{PBGclass}. Recall that the representation $\rho$ induces an atlas of PBG-Lie group bundle charts $\rho(\sigma_{i}(u),h)$ for $F \rightarrow P(M,G)$. Thus every element of $F$, say $\lambda \in F_{u}$, can be represented as $\rho(\sigma_{i}(u),h)$ for any $i \in I$ with $u \in P_{i}$. Define $\iota : F \rightarrow \Upsilon$ by mapping $\rho(\sigma_{i}(u),h) \in F_{u}$ to $\langle i,(u,h,u) \rangle$. It is trivial to check that $\iota$ is well defined, and an isomorphism of PBG-Lie group bundles over $P(M,G)$ onto $I\Upsilon$. Thus we have the extension of PBG-groupoids
\[
F \stackrel{\iota}{\inj} \Upsilon \stackrel{(\beta,\alpha)}{\surj} P \times P.
\]
Define $\mu : \Upsilon \rightarrow \Omega$ by $\langle i, (u,h,v) \rangle \mapsto \sigma_{i}(u)\da(h)\sigma_{i}(v)^{-1}$. Again one checks that $\mu$ is well defined, a surjective submersion and a morphism of PBG-groupoids over $P(M,G)$.

To see that the diagram
\begin{diagram}
F & \rto^{\iota} & \Upsilon \\
\dTeXonto<{\da} &  & \dTeXonto<{\mu} \\
Im\da & \rto & \Omega
\end{diagram}
commutes, recall that $\da(\rho(\xi,\lambda)) = \xi \da(\lambda) \xi^{-1}$ for $\xi \in \Omega, \lambda \in F_{\alpha\xi}$. Taking $\xi = \sigma_{i}(u)$ and $\lambda = h \in H = F_{u_{0}}$, this gives
\[
\da(\rho(\sigma_{i}(u),h)) = (\mu \circ \iota)(\rho(\sigma_{i}(u),h)),
\]
as required. 

It remains to verify that the action of $\Omega$ on $F$ induced by the diagram 
\begin{diagram}
\ker\da & \requal & \ker\da &  &  \\
\dTeXembed &  & \dTeXembed &  & \\
F & \rTeXembed^{\iota} & \Upsilon & \rTeXonto^{(\beta,\alpha)} & P \times P \\
\dTeXonto<{\da} &  & \dTeXonto<{\mu} &  & \dequal \\
Im\da & \rTeXembed & \Omega & \rTeXonto^{(\beta,\alpha)} & P \times P
\end{diagram}
coincides with the given $\rho$. Take $\omega \in \Upsilon$, say $\omega = \langle j, (u,h,v) \rangle$, and $\lambda \in F_{\alpha\omega}$, say $\lambda = \rho(\sigma_{i'}(u),h')$; it is no loss of generality to assume that $j = i'$. Now $\omega \iota(\lambda) \omega^{-1} = \langle j, (u,hh'h^{-1},v) \rangle$, by the definition of $\iota$ and the multiplication in $\Upsilon$. On the other hand, $\mu(\omega)$ is equal to $\sigma_{j}(u)\da(h)\sigma_{j}(v)^{-1}$ and so 
\[
\rho(\mu(\omega),\lambda) = \rho(\sigma_{j}(u)\da(h),h') = \rho(\sigma_{j}(u),hh'h^{-1}).
\]
So $\omega\iota(\lambda)\omega^{-1} = \iota(\rho(\mu(\omega),\lambda))$, as required. This completes the proof that $(\Upsilon,\mu)$ is an operator extension of the pair crossed module of PBG-groupoids $\langle F,\da,\Omega,\rho \rangle$. The converse is a trivial verification. \boom

The element $e \in \check{H}^{2}_{G}(P \times P,ZF)$ is the {\em obstruction} associated with the pair crossed module of PBG-groupoids $\langle F, \da,\Omega,\rho \rangle$. Following the notation of \cite{Mackenzie:1989}, we denote it by $\Obs\langle F, \da,\Omega,\rho \rangle$. The following theorem is an immediate consequence of the previous considerations.
\begin{thm}
Let $\Omega \gpd M$ be a transitive Lie groupoid and 
\begin{eqnarray}\label{thm:integrobstr}
K \inj A \surj A\Omega
\end{eqnarray}
be an extension of Lie algebroids over the manifold $M$. Choose a baspoint in $M$ and let $P(M,G,p)$ be the principal bundle corresponding to $\Omega$. If $F \rightarrow P(M,G)$ is the PBG-Lie group bundle integrating the PBG-Lie algebra bundle $K \ract p$, then the integrability obstruction of the extension (\ref{thm:integrobstr}) is the obstruction $e \in \check{H}^{2}_{G}(P \times P, ZF)$ associated with the pair crossed module of PBG-groupoids associated to (\ref{thm:integrobstr}).
\end{thm}

\section{Classification of operator extensions for coupling pair crossed modules of PBG-groupoids}

Suppose given a coupling pair crossed module of PBG-groupoids $\langle F, \da, \Omega, \rho \rangle$ over the principal bundle $P(M,G)$. Recall that {\em coupling} means $\ker\da = ZF$ whereas {\em pair} means that the cokernel $\frac{\Omega}{\ker\da} = P \times P$. In this section we show that if its obstruction cocycle vanishes then its operator extensions are classified by $\check{H}^{1}_{G}(P, ZH)$, where $H$ is the fiber type of the Lie group bundle $F$. This cohomology, defined on $P$ instead of $P \times P$ is the isometablic cohomology given in \cite{Androulidakis:2004CELG}, and we start with a brief recollection of it.

Consider a PBG-groupoid $\Xi \gpd P(M,G)$. Choose an atlas $\{ P_{i} \equiv U_{i} \times G \}_{i \in I}$ for the principal bundle $P(M,G)$, where $\{ U_{i} \}_{i \in I}$ is a simple open cover of $M$, and a family of local flat isometablic connections $TP_{i} \rightarrow A\Xi_{P_{i}}$. As we discussed in the beginning of section 6, this data gives rise to sections $\sigma_{i} : P_{i} \rightarrow \Xi_{u_{0}}$ of the PBG-groupoid, which are isometablic in the sense
\[
\sigma_{i}(ug) = [\sigma_{i}(u)g]\cdot (\xi_{i}^{-1}g) \cdot \sigma_{i}(u_{i}g).
\]
An alternative classification of PBG-groupoids, given in \cite{Androulidakis:2004CELG}, is by the transition functions $\{ s_{ij} : P_{ij} \rightarrow \Xi_{u_{0}}^{u_{0}} \}_{i,j \in I}$ of these sections. The isometablicity of these functions is expressed by
\[
s_{ij}(ug) = \phi_{ij}(g)(s_{ij}(u)).
\]
Here, the $\phi_{ij}$s are the actions of $G$ on $\Xi_{u_{0}}^{u_{0}}$ defined by
\[
\phi_{ij}(g)(h) = \sigma_{i}(u_{i}g)^{-1} \cdot (\xi_{i}g) \cdot (hg) \cdot (\xi_{j}g)^{-1} \cdot \sigma_{j}(u_{j}g).
\]
Note that the $\phi_{ij}$s are just actions, {\em not} representations of $G$ on $\Xi_{u_{0}}^{u_{0}}$. That is because for every $g \in G$, the map $\phi_{ij}(g) : \Xi_{u_{0}}^{u_{0}} \rightarrow \Xi_{u_{0}}^{u_{0}}$ does not preserve the multiplication on $\Xi_{u_{0}}^{u_{0}}$. Instead, it is a straightforward calculation that for $P_{ijk} \neq \emptyset$ they satisfy
\[
\phi_{ij}(g)(h_{1}h_{2}) = \phi_{ik}(g)(h_{1})\phi_{kj}(g)(h_{2}).
\]
This property was called a {\em cocycle morphism} in \cite{Androulidakis:2004CELG}, and it was shown that such data (that is to say cocycle morphisms $\phi = \{\phi_{ij} : G \times H \rightarrow H \}_{i,j \in I}$, together with a $\phi$-isometablic cocycle $\{ s_{ij} : P_{ij} \rightarrow H \}_{i,j \in I}$, where $H$ is a Lie group) classifies PBG-groupoids. 



Now the definition of isometablic \v{C}ech cohomology with respect to the family of actions $\phi_{ij}$ was given in \cite[\S VII]{Androulidakis:2004CELG}, and it was shown that $\check{H}^{1}_{G}(P,H)$ classifies those PBG-groupoids over the principal bundle $P(M,G)$ such that the fiber of the adjoint bundle is the Lie group $H$.

Let us make a fresh start now, considering a coupling pair crossed modules of PBG-groupoids $\langle F, \da, \Omega, \rho \rangle$ over the principal bundle $P(M,G)$. Let $H$ denote the fiber type of the PBG-Lie group bundle $F$, and suppose that $\Obs\langle F, \da,\Omega,\rho \rangle = 0$. Let $\Opext\langle F, \da,\Omega,\rho \rangle$ denote the set of equivalence classes of operator extensions. We define an action of $\check{H}^{1}_{G}(P \times P,ZH)$ on $\Opext\langle F, \da,\Omega,\rho \rangle$ in the following way:

Consider an operator PBG-groupoid $(F \stackrel{\iota}{\inj} \Upsilon \stackrel{(\beta,\alpha)}{\surj} P \times P,\mu)$ for $\langle F, \da, \Omega, \rho \rangle$, and an element $f \in H^{1}_{G}(P \times P,ZH)$. Note that $f : P_{ij} \times P_{ij} \rightarrow ZH$ is isometablic in the sense
\[
f_{ij}(ug,vg) = \phi_{ij}(g)(f_{ij}(u,v))
\]
for all $g \in G$. Therefore, if $\Hat{s}_{ij}$ are the transition functions of the PBG-groupoid $\Upsilon$, arising from an isometablic section-atlas $\Hat{\sigma}_{i} : P_{i} \rightarrow \Upsilon_{u_{0}}$, the maps $\Hat{s}_{ij}f_{ij} : P_{ij} \rightarrow H$ satisfy the cocycle equation and $[s_{ij}f_{ij}](ug) = \phi_{ij}(g)(s_{ij}f_{ij}(u))$. Moreover, $\da \circ (\Hat{s}_{ij}f_{ij}) = \da \circ \Hat{s}_{ij}$.
\begin{prop}
The PBG-groupoid $\Upsilon^{f} \gpd P(M,G)$ constructed from the $\Hat{s}_{ij}f_{ij}$s, is an operator PBG-groupoid for $(F \stackrel{\iota}{\inj} \Upsilon \stackrel{(\beta,\alpha)}{\surj} P \times P,\mu)$.
\end{prop}

\pf Same as \cite[3.4]{Mackenzie:1989}. \boom

The proof that this action is well defined is exactly the same as in \cite[\S3]{Mackenzie:1989}, taking into account the isometablicity considerations of section 6 in the present paper. Moreover, applying these considerations to the proof of \cite[3.5]{Mackenzie:1989}, we get the following classification of operator extensions for a coupling pair crossed module of PBG-groupoids:
\begin{thm}
The above action of $\check{H}^{1}_{G}(P \times P,ZH)$ on the set of operator extensions of a coupling pair crossed module of PBG-groupoids $\langle F, \da, \Omega, \rho \rangle$ is free and transitive.
\end{thm}


\nocite{Brown-Mackenzie:1992DLG} \nocite{Androulidakis:solo1} \nocite{Moerdijk:classification}

\bibliographystyle{plain}


\end{document}